\documentclass[11pt]{amsart}
\usepackage{amssymb}
\usepackage{amsfonts}
\usepackage{amscd}
\usepackage{amsmath}
\usepackage{mathrsfs}
\usepackage{curves}
\usepackage{epsfig}
\usepackage[pass]{geometry}
\usepackage{graphicx}
\usepackage{verbatim}
\usepackage{latexsym}
\usepackage{eucal}
\usepackage{bbm}
\usepackage{tikz}
\usetikzlibrary{matrix}

\numberwithin{equation}{section}

\newtheorem{theorem}{Theorem}[section]
\newtheorem{lemma}[theorem]{Lemma}
\newtheorem{prop}[theorem]{Proposition}
\newtheorem{cor}[theorem]{Corollary}

\newtheorem{prob}[theorem]{Problem}

\def \mcc {{\mathscr C}}

\def \mck {{\mathscr K}}
\def \mcl  {\mathcal{L}}
\def \mcm {{\mathscr M}}

\def \mcs {{\mathscr S}}

\def \mcu {{\mathscr U}}
\def \mcv {{\mathscr V}}

\def \mbn {{\mathbb N}}
\def \mbr {{\mathbb R}}
\def \mbs {{\mathbb S}}

\def \id {\operatorname{Id}}
\def \comp {\operatorname{comp}}
\def \loc {\text{loc}}

\def \WF {\text{WF}}

\def \supp {\text{supp }}
\def \defeq {\stackrel{\operatorname{def}}{=}}
\def \beqq {\begin{equation}}
\def \eeqq {\end{equation}}
\def \bpf {\begin{proof}}
\def \epf {\end{proof}}
\def \beq {\begin{equation*}}
\def \eeq {\end{equation*}}

\def \eps {\epsilon}   
\def \la {\lambda}

\def \lap {\Delta}
\def \p {\partial}
\def \ha {\frac{1}{2}}
\def \tilde {\widetilde}

 \def \mcm {\mathscr{M}}

 \def \intm  {M^\circ}

\begin{document}
\title[]{On the light ray transform of wave equation solutions}
\author{Andr\'as Vasy}
\address{Andr\'as Vasy
\newline
\indent Department of Mathematics, Stanford University}
\email{andras@stanford.edu}
\author{Yiran Wang}
\address{Yiran Wang 
\newline
\indent Department of Mathematics, Emory University}
\email{yiran.wang@emory.edu}
\begin{abstract} 
We study the light ray transform on Minkowski space-time and its small metric perturbations acting on scalar functions which are solutions to wave equations.  We show that the light ray transform uniquely determines the function in a stable way.  The problem is of particular interest because of its connection to inverse problems of the Sachs-Wolfe effect in cosmology.
\end{abstract}
\maketitle

\section{Introduction}\label{}
\begin{figure}[htbp]
 \centering
\includegraphics[scale = 0.55]{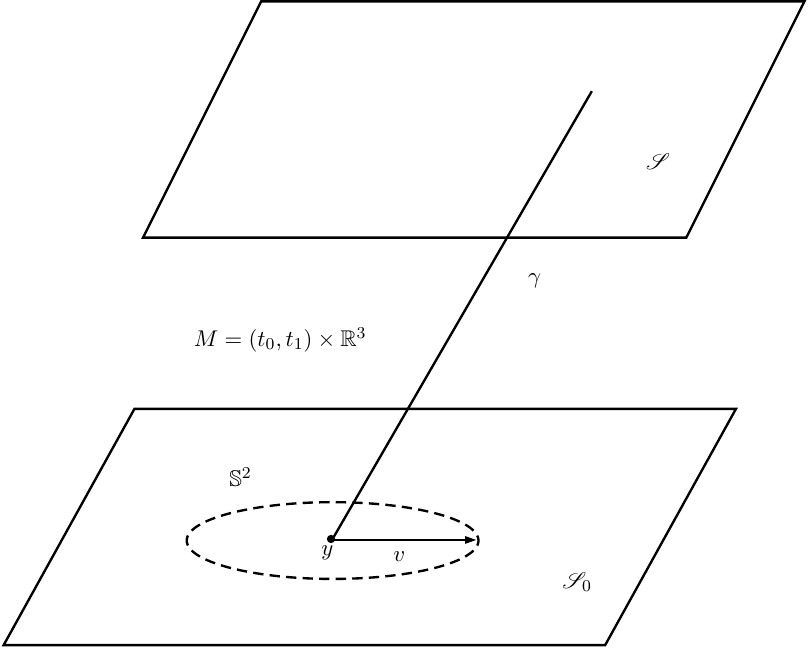}
\caption{The setup of the problem for the Minkowski space-time.}
\label{fig-setup}
\end{figure}

Let $M = [t_0, t_1]\times \mbr^3$ and $(t, x), t\in [t_0, t_1], x\in \mbr^3$ be the local coordinates. Let  $g_M = -dt^2 + dx^2$ be the Minkowski metric on $M$. Consider the Lorentzian manifold $(M, g_M)$. We denote the interior by 
  $\intm = (t_0, t_1)\times \mbr^3$ and  the boundaries by $\mcs_0 = \{t_0\}\times \mbr^3$ and $\mcs = \{t_1\}\times \mbr^3.$ See Figure \ref{fig-setup}.   

Consider light-like geodesics on $(M, g_M)$ which are straight lines. We parametrize the set of light rays $\mcc$ as follows: let $x_0\in \mcs_0$ and $v \in \mbs^2$ the unit sphere in $\mbr^3$. Then a light ray from $x_0$ in direction $(1, v)$ is $\gamma(\tau) = (t_0, x_0) + \tau(1, v),  \tau \in [0, t_1 - t_0].$ See Figure \ref{fig-setup}. In particular, we can identify $\mcc = \mbr^3\times \mbs^2.$ 
The light ray transform for scalar functions on $(M, g_M)$ is defined by 
\beqq\label{eq-lray}
X_M(f)(\gamma) = \int_{0}^{t_1-t_0}f(\gamma(\tau)) d\tau, \ \ f\in C_0^\infty(M).
\eeqq
Of course, one can regard $X_M$ as the restriction of the light ray transform $X_{\mbr^4}$ of the Minkowski spacetime $(\mbr^4, g_M)$ acting on functions supported in $M$. However, it is perhaps better to think of $X_M$ as the compact version of the transform, which is similar to the geodesic ray transform on a compact Riemannian manifold with boundary, see for instance \cite{Sha}. 

In this work, we study $X_M$ acting on scalar functions which are solutions to the Cauchy problem of wave equations on $M$. 
Let $c > 0$ be a constant. Denote $\square_c = \p_t^2 + c^2 \lap $ where $\lap$ is the positive Laplacian on $\mbr^3,$ namely $\lap = \sum_{i = 1}^3 D_{x_i}^2, D_{x_i} = -\sqrt{-1} \frac{\p}{\p x_i}.$ Here, $c$ is the wave speed. On  $(M, g_M)$,  $c = 1$ is the speed of light, and $\square_c$ is the d'Alembert operator. Consider the Cauchy problem
\beqq\label{eq-cons}
\begin{gathered}
\square_c f = 0 \ \ \text{ on } M\\
f = f_1, \ \ \p_t f = f_2, \ \ \text{ on } \mcs_0.
\end{gathered}
\eeqq
 The problem we address in this paper is the determination of $f$ or equivalently $f_1, f_2$ from $X_M(f)$ with the constraint \eqref{eq-cons}.  Let $\mathcal{N}^s \defeq  H_{\comp}^{s+1}(\mcs_0) \times H_{\comp}^s(\mcs_0)$. 
 Our main result is 
 \begin{theorem}\label{thm-main1}
Suppose $0 < c \leq 1$ is constant. Assume that $(f_1, f_2)\in \mathcal{N}^s, s\geq 0$, and $f_1, f_2$ are supported in a compact set $\mck$ of $\mcs_0$. Then $X_M f$ uniquely determines $f$  and $f_1, f_2$ which satisfy \eqref{eq-cons}. Moreover, there exists   $C> 0$ such that 
\beq
\|(f_1, f_2)\|_{\mathcal{N}^s} \leq C  \|X_M  f\|_{H^{s+ 2}(\mcc)} \text{ and } \|f\|_{H^{s+1}(M)} \leq C  \|X_M   f\|_{H^{s+ 2}(\mcc)} 
\eeq
where  $\mcc$ is the set of light rays on $M$.
 \end{theorem}
 
We will prove stronger versions of the theorem including lower order  terms in the wave equation in Theorem \ref{thm-main2} in Section \ref{sec-full}.  However, for ease of presentation, we use the
 standard wave equation on Minkowski spacetime throughout the paper until
 the final sections where the necessary changes are indicated.  
 
Next, we consider the generalization of Theorem \ref{thm-main1} corresponding to $c = 1$. We remark that it is not difficult to formulate the result corresponding to $c<1$ although we do not discuss it. We consider metric perturbations $g_\delta = g_M + h$ where $h$  satisfies assumptions (A1), (A2) in Section \ref{sec-pert}, which says that $h$ is a suitably smooth small perturbation of the Minkowski spacetime. In this case, light rays may not be straight lines. Let $X_\delta$ be the light ray transform on $(M, g_\delta)$ see \eqref{eq-xraypert}. Let $\square_{g_\delta}$ be the d'Alembert operator on $(M, g_\delta)$. Consider the Cauchy problem
\beqq\label{eq-cons-1}
\begin{gathered}
\square_{g_\delta} f = 0 \ \ \text{ on } \intm\\
f = f_1, \ \ \p_t f = f_2, \ \ \text{ on } \mcs_0.
\end{gathered}
 \eeqq 
 Our result is  
 \begin{theorem}\label{thm-main-1}
Consider $(M, g_\delta)$ satisfying   assumptions (A1), (A2) to be stated in  Section \ref{sec-pert}. Assume that $(f_1, f_2)\in \mathcal{N}^s, s \geq 0$, and $f_1, f_2$ are supported in a compact set $\mck$ of $\mcs_0$. For $\delta \geq 0$ sufficiently small, $X_\delta f$ uniquely determines $f$ and $f_1, f_2$ which satisfy \eqref{eq-cons-1}. Moreover, there exists   $C> 0$ such that 
\beq
\|(f_1, f_2)\|_{\mathcal{N}^s} \leq C  \|X_\delta   f\|_{H^{s+ 2}(\mcc_\delta)} \text{ and } \|f\|_{H^{s+1}(M)} \leq C  \|X_\delta  f\|_{H^{s+ 2}(\mcc_\delta)} 
\eeq
where  $\mcc_\delta$ is the set of light rays on $(M, g_\delta)$, see Section \ref{sec-pert}.
 \end{theorem}

Our motivation for this setup of the light ray transform comes from some inverse problems in cosmology. We are particularly interested in the determination of gravitational perturbations such as primordial gravitational waves from the anisotropies of the Cosmic Microwave Background (CMB), see for example \cite{KDS, Do, Dur}.  Sachs and Wolfe in their  1967 paper   \cite{SW} discovered the connection of the CMB anisotropy and the light ray transform of the gravitational perturbations, now called the Sachs-Wolfe effects. We discuss the background in Section \ref{sec-sw} and \ref{sec-wave}. Physically, $c< 1$ and $c  =1$ in Theorem \ref{thm-main1} correspond  to different Universe models driven by hydrodynamical perturbations and scalar field perturbations, respectively. Moreover, Theorem \ref{thm-main-1} covers some cases of variable wave speeds.  
 
The reason that we are able to get a stable determination is the restriction of singularities of $f$. In general, it is known that time-like singularities in $f$, namely all $(z, \zeta) \in T^*M$ in the wave front set $\WF(f)$ of $f$ with $\zeta$ time-like, are lost after taking the light ray transform,  
although  the light ray transform $X_M$ is injective on $C_0^\infty(M)$. In particular, we do not expect Theorem \ref{thm-main1} to hold for $c>1$. There is a fundamental difference in our treatment between the $c < 1$ and $c = 1$ cases. The former requires a good understanding of the normal operator $X_M^*X_M$ which was considered in \cite{LOSU} and further generalized in \cite{LOSU1}, while the latter relies on a thorough analysis of the operator $X_M E$ where $E$ is the fundamental solution or parametrix for the Cauchy problem.  

The paper is organized as follows. In Section \ref{sec-sw} and \ref{sec-wave}, we discuss the (integrated) Sachs-Wolfe effects and explain how the inverse problem is related to our theorems. In Section \ref{sec-light}, we review some properties of the light ray transform. Then we consider the Cauchy problem in Section \ref{sec-cauchy}. In Section \ref{sec-local1} and \ref{sec-local2}, we construct the microlocal parametrix for the light ray transform with the wave constraint for $c< 1$ and $c = 1$ respectively. We prove  Theorem \ref{thm-main1} and the version including lower order terms in the wave equation in Section \ref{sec-full}. Finally, we address the small metric perturbations of Minkowski space-time in Section \ref{sec-pert}.

\section{The integrated Sachs-Wolfe effect}\label{sec-sw}
Consider the flat Friedman-Lema\^ite-Robertson-Walker (FLRW) model for the cosmos:
\beq
\mcm = (0, \infty)\times \mbr^3, \ \ g_0 = dt^2 - a^2(t) \delta_{ij} dx^i dx^j
\eeq
where $(t, x), t \in (0, \infty), x\in \mbr^3$ are coordinates and $\delta_{ij} = 1$ if $i = j$ and otherwise $0$. Here, the signature of $g_0$ is $(+, -, -, -)$ because we will refer to some results in \cite{MFB} later.  The factor $a(t)$ is assumed to be positive and smooth in $t$. It represents the rate of expansion of the Universe.  

We assume that the actual cosmos is a metric perturbation $
g = g_0 + \delta g$  on $\mcm$ 
where $\delta g$ is a small perturbation compared to $g_0.$  Here, we follow  the convention of \cite{MFB} that $\delta A$ denotes the perturbation of quantity $A$ (not $\delta$ times $A$).  We introduce the conformal time $s$ such that $ds = a^{-1}(t)dt$. Then we get 
\beq
g_0 = a^2(s) (ds^2 - \delta_{ij} dx^i dx^j) = a^2(s) g_{M}
\eeq
where $g_{M}$ is the Minkowski metric on $\mcm = (0, \infty)$ and we used $a(s)$ to denote $a(t(s))$. We write  $
g = a^2(s) (g_{M} + \delta g)
$ 
where $\delta g$ denotes the corresponding perturbation in conformal time.  In the literature, the metric perturbations are classified to scalar, vector and tensor type. We consider the scalar type perturbations. In the longitudinal gauge, also called the  conformal Newtonian gauge, the metric $g$ is of the form
\beqq\label{eq-confg}
g = a^2(s)[ (1+ 2\Phi) ds^2 - (1-2\Psi) dx^2]
\eeqq
 see \cite[Section 2]{MFB}. Here, $\Phi, \Psi$ are scalar functions on $M.$ We remark that there is a gauge invariant formulation of cosmological perturbations. However, in the longitudinal gauge, the gauge invariant variables are equal to $\Phi, \Psi$, see \cite{MFB}. In this work, we fix the gauge and work with $\Phi, \Psi$ for simplicity.  
 
Consider the  Cosmic Microwave Background (CMB) measurement. Our main references are \cite{Do, Dur, SW}. Let $\mcs_0 = \{s_0\}\times \mbr^3$ be the surface of last scattering. This is the moment after which photons stopped interaction and started to travel freely in $\mcm.$ Let $\mcs = \{s_1\}\times \mbr^3$ be the surface where we make observation of the photons. Let $\gamma(\tau)$ be a light ray from $\mcs_0$ to $\mcs$. It represents the trajectory of photons in $\mcm.$ Explicitly, we have 
\beq
\gamma(\tau) = (s_0, x_0) + \tau (1, v), \ \ (s_0, x_0) \in \mcs_0, v\in \mbs^2, \tau \in [0, s_1-s_0].
\eeq
Then we consider the photon energies observed  at $\mcs_0, \mcs$ denoted by  
$
E_0 = g_0(\dot \gamma(s_0), \p_s), E = g_0(\dot \gamma(s_1), \p_s). 
$
Here, the observer is represented by the flow of the vector field $\p_s.$ The redshift $z$ is defined by 
\beq
1 + z = E/E_0.
\eeq
In \cite{SW}, Sachs and Wolfe derived that to the first order linearization, $1+ z$ is represented by a light ray transform of the metric perturbations, see \cite[equation (39)]{SW}. In cosmological literatures, one often connects this to the CMB temperature anisotropies. Let $T$ be the temperature observed at $\mcs$ in the isotropic background $g_0$. Let $\delta T$ be the temperature fluctuation from the isotropic background. One can compute $\delta T/T$ in terms of the energies $E_0, E$. One component of $\delta T/T$ is the integrated Sachs-Wolfe (ISW) effects 
\beqq\label{eq-eqsw}
\begin{gathered}
(\frac{\delta T}{T})^{ISW}(\gamma) = \int_{0}^{s_1-s_0}(\p_s \Phi(\gamma(\tau)) + \p_s \Psi(\gamma(\tau)) d\tau 
 = X_M (\p_s \Phi + \p_s \Psi)(\gamma)
\end{gathered}
\eeqq
 see \cite[Section 2.5]{Dur}. Note that this quantity depends on the light ray $\gamma$ which indicates the anisotropy. 
  We remark that another component of $\delta T/T$ is the ordinary Sachs-Wolfe effect (OSW) which only involves $\Phi, \Psi$ at $\mcs_0$. The integrated Sachs-Wolfe effect can be extracted from the CMB and other astrophysical data, see for example \cite{MD}.  
 
The inverse Sachs-Wolfe problem we study is to determine $\Phi, \Psi$ on $M$ from $(\delta T/T)^{ISW}$, which in particular includes  the initial value of $\Phi, \Psi$ on $\mcs_0$. 
Before we proceed, we observe that there are natural obstructions to the unique determination from \eqref{eq-eqsw}. If $\Phi + \Psi$ is a constant, then the integrated Sachs-Wolfe effect is always zero. So the goal is to determine $\Phi, \Psi$ up to such natural obstructions.

\section{Dynamical equations for perturbations}\label{sec-wave}
For the Sachs-Wolfe problem, we should take into account that $g$ satisfies the Einstein equations with certain source fields and initial perturbations at $\mcs_0$ from $g_0$. On the linearization level, this puts the perturbation $\delta g$  under some wave equation constraint as we discuss in this section. The derivations of the equations for the perturbation take some amount of work and they are mostly done in the literature, see for example \cite[Section 5.1]{Do} and \cite{Dur}. We follow the presentation and the notations in \cite[Section 4-6]{MFB} closely.  Instead of the gauge invariant approach, we choose to work in the longitudinal gauge for simplicity. It is not hard to transform back and forth and our analysis works for the gauge invariant formulation as well.  

Let $R^\mu_{\ \ \nu}$ be the Ricci curvature tensor and $R$ the scalar curvature on $(\mcm, g)$ (in conformal time).  Let $T^\mu_{ \ \ \nu}$ denote the stress-energy tensor of certain source fields.  The Einstein equations are 
\beq
G^\mu_{\ \ \nu} = 8\pi G T^\mu_{\ \ \nu}, \ \ G^\mu_{\ \ \nu} = R^\mu_{\ \ \nu} - \ha \delta^\mu_{\ \ \nu} R
\eeq
where $G$ is Newton's gravitational constant. 
We assume that $T^\mu_{\ \ \nu} = {}^{(0)}T^\mu_{\ \ \nu} + \delta T^\mu_{\ \ \nu}$ where ${}^{(0)}T$ denotes the stress-energy tensor of the background field and $\delta T$ denotes the perturbation. We also have  $g = a^2 (g_M + \delta g)$. Then we can write $G^\mu_{\ \ \nu} = {}^{(0)}G^\mu_{\ \ \nu} + \delta G^\mu_{\ \ \nu} + \cdots. $ From the asymptotic expansion, one finds that the Einstein tensor for the background metric $g_M$ are 
\beqq\label{eq-sderi}
\begin{gathered}
{}^{(0)}G_0^{\ \ 0} = 3a^{-2} H^2, \ \ {}^{(0)}G_{\ \ i}^0 = 0,  \ \ 
{}^{(0)}G_{\ \ j}^i = a^{-2}(2H' + H^2) \delta^i_{\ \ j},  
\end{gathered}
\eeqq
where $i, j = 1, 2, 3, $ $H(s) = \p_s a(s)/a(s)$, see \cite[equation (4.2)]{MFB}. Here, $H' = \p_s H$ denotes the derivative in the conformal time variable. We emphasize that we work with a flat Universe and we get the equation ${}^{(0)}G^\mu_{\ \ \nu} = 8\pi G {}^{(0)}T^\mu_{\ \ \nu}$. 

For the first order perturbation term, we get $\delta G^\mu_{\ \ \nu} = 8\pi G \delta T^\mu_{\ \ \nu}$. After lengthy calculations, one obtains (see \cite[equation (4.15)]{MFB}) the following equations for $\Phi, \Psi$ 
\beqq\label{eq-eqeinlin}
\begin{gathered}
-3H (H \Phi + \Psi') + \lap \Psi  = 4\pi G a^2 \delta T^0_{\ \ 0}\\
\p_i (H \Phi + \Psi')  = 4\pi G a^2 \delta T^0_{\ \ i}\\
[(2 H' + H^2)\Phi + H \Phi' + \Psi'' + 2H \Psi' + \ha \lap (\Phi - \Psi) ] \delta^i_{\ \ j} - \ha \delta^{ik}(\Phi - \Psi)_{|kj} = -4\pi Ga^2 \delta T^i_{\ \ j},
\end{gathered}
\eeqq
where $i, j = 1, 2, 3$, $\p_i$ denotes the $i$th component of the covariant derivative with respect to the background metric $g_M$, $\lap$ denotes the standard Laplacian on $\mbr^3$, and as in \eqref{eq-sderi}, prime denotes $\p_s$ derivative.  

Now we need to specify the source field. We consider two important examples: the perfect fluid and the scalar field. 
We first consider Universe dominated by perfect fluid sources. Let $u$ be the four fluid velocity of a fluid source. The stress-energy tensor for a perfect fluid is 
\beq
T^\alpha_{\ \ \beta} = (\eps + p) u^\alpha u_\beta - p\delta^\alpha_{\ \ \beta}
\eeq
see \cite[equation (5.2)]{MFB}, Here, $\eps$ is the energy density and $p$ is the pressure of the fluid. We assume that $\eps = \eps_0 + \delta \eps, p = p_0 +  \delta p$ where $0$ denotes the quantity for the background and $\delta$ denotes the perturbations. For fluid source, from \eqref{eq-eqeinlin} one deduces that the perturbations $\Phi = \Psi$. In the case of adiabatic perturbations, $\Phi$ satisfies the following equation, called Bardeen's equation
\beqq\label{eq-bardeen}
\Phi'' + 3H(1 + c_s^2)\Phi' - c_s^2 \lap \Phi + [2H' + (1 + 3c_s^2)H^2 ]\Phi = 0,
\eeqq
see \cite[equation (5.22)]{MFB}. 
In general, the right hand side of the equation is a non-zero term related to the entropy perturbations. The fluid velocity $u$ also satisfies a wave equation with speed $c_s$, see \cite[equation (5.25)]{MFB}. Here, $c_s < 1$ is the speed of sound. Prescribing Cauchy data of $\Phi$ at $\mcs_0$, one can solve the Cauchy problem of \eqref{eq-bardeen} to get $\Phi$ in $\mcm$. We formulate the inverse Sachs-Wolfe problem in this case as 
\begin{prob}
Determining $\Phi$ from \eqref{eq-eqsw} where $\Phi$ satisfies the Cauchy problem of \eqref{eq-bardeen}.
\end{prob}
Commuting equation \eqref{eq-bardeen} with $\p_s$, we see that $\p_s\Phi$ also satisfies a wave equation. Hence, we arrived at the model problem we proposed in the introduction.

Next, let's consider Universe governed by a scalar field $\phi$. The stress energy tensor is  
\beq
T^\mu_{\ \ \nu} = \nabla^\mu \phi \nabla_\nu \phi - [\ha \nabla^\alpha\phi \nabla_\alpha\phi - V(\phi)]\delta^\mu_{\ \ \nu}
\eeq
see \cite[equation (6.2)]{MFB}. Here, 
 $V$ is the potential function for the scalar field $\phi$. The field itself satisfies the Klein-Gordon equation $
\square \phi + \p_\phi V(\phi) = 0.$ 
Now assume that $\phi = \phi_0 + \delta \phi$ where $\phi_0$ is the scalar field which drives the background model and $\delta \phi$ denotes the perturbation. Then we can split $T^\mu_{\ \ \nu} = {}^{(0)}T^\mu_{\ \ \nu} + \delta T^{\mu}_{\ \ \nu}$. Again, one finds that $\Phi = \Psi$ and it satisfies the equation 
\beqq\label{eq-eqscalar}
\Phi'' + 2(H - \phi_0''/\phi_0') \Phi' - \lap \Phi + 2(H' - H \phi_0''/\phi_0)\Phi = 0
\eeqq
see \cite[equation (6.48)]{MFB}. This is a damped wave equation with wave speed $c = 1$. We can formulate the inverse Sachs-Wolfe problem in this case as
\begin{prob}
Determining $\Phi$ from \eqref{eq-eqsw} in which $\Phi$ satisfies the Cauchy problem of \eqref{eq-eqscalar}. 
\end{prob}
Again, we arrived at the model problem in the introduction with $c = 1. $ We do not need it but record that the scalar field perturbation also satisfies a wave equation, 
see \cite[equation (6.47)]{MFB}. 

Applying our main result of the paper, in particular Theorem \ref{thm-main2} which allows  lower order terms in the wave equation, we obtain the following result. 
\begin{cor}
For the inverse Sachs-Wolfe effect Problems 3.1 and 3.2, one can uniquely determine $\Phi$ in $\mcm$ (and the initial conditions at $\mcs_0$) in the longitudinal gauge up to a constant in a stable way. 
\end{cor}

\section{The light ray transform on functions}\label{sec-light}
We recall some facts about the light ray transform on scalar functions. Consider the Lorentzian manifold $(M, g_M)$ and hereafter we change the signature of $g_M$ to $(-, +, +, +).$ For $ (t, x)\in \intm, t\in (t_0, t_1), x\in \mbr^3$, we use $\Xi = (\tau, \xi), \tau\in \mbr, \xi\in \mbr^3$ for the coordinate in  $T_{(t, x)}\intm$ so that tangent vectors are represented by $\tau\p_t + \sum_{j = 1}^3 \xi^j\p_{x^j}$. We divide the tangent vectors in $T_{(t, x)}\intm$ into time-like vectors $\Omega^-_{(t, x)}  \intm = \{\Xi \in \mbr^4:  g_M(\Xi, \Xi) = -\tau^2+ |\xi|^2 < 0 \}$, space-like vectors $\Omega^+_{(t, x)} \intm = \{\Xi \in \mbr^4: g_M(\Xi, \Xi) > 0 \}$ and light-like vectors $L_{(t, x)} \intm = \{\Xi \in \mbr^4: g_M(\Xi, \Xi) = 0 \}$. We denote the corresponding fiber bundles by $\Omega^-\intm, \Omega^+\intm, L\intm.$ The cotangent vectors can be classified similarly using the dual metric $g_M^*$ on $T^*\intm.$ The corresponding bundles are denoted by $\Omega^{*,-}\intm, \Omega^{*, +}\intm, L^*\intm.$

From now on, {\em without loss of generality, we take $t_0 = 0$ in $M$}, which amounts to a translation in the $t$ variable. 
Let $\mcc$ be the set of light rays on $(M, g_M)$. As $M$ has a global coordinate system, we can  parametrize $\mcc$ as follows. Let $y\in   \mbr^3, v \in \mbs^2 \defeq \{z \in \mbr^3: |z| = 1\}$ with $|\cdot|$ the Euclidean norm. We denote $\theta = (1, v)$ so that $\theta$ is a (future pointing) light-like vector. Then all the light rays are given by $\gamma_{y, v}(\tau) =   (\tau, y + \tau v),  \tau \in (0, t_1), (y, v)\in \mbr^3\times \mbs^2$. Thus, we can identify $\mcc  = \mbr^3\times \mbs^2.$ 
For $f\in C_0^\infty(\intm)$ and $ y\in \mbr^3, v\in \mbs^2$, we have
\beqq\label{eqlit}
\begin{gathered}
X_{M}f(y, v) = \int_0^{t_1} f(\tau, y + \tau v)d\tau 
 = (2\pi)^{-3}\int_{\mbr^3}\int_{\mbr^3}  \int_{0}^{t_1} e^{i ( (y - x)\cdot \eta + t v\cdot \eta)} f(t, x)  dt dx  d\eta.
\end{gathered}
\eeqq
The Schwartz kernel of $X_{M}$ is $\delta_Z$ the delta distribution on $\mcc \times \intm$ supported on the point-line relation $Z$ defined by  
\beq
\begin{gathered}
Z = \{(\gamma, q)\in \mcc\times \intm: q\in \gamma\} 
= \{(y, v, (t, x)) \in \mbr^3 \times \mbs^2 \times \intm:  x = y+ t  v\}.
\end{gathered}
\eeq
We know (see e.g.\ \cite{LOSU}) that $X_{M}$ is an Fourier integral operator of order $-3/4$ associated with the canonical relation $(N^*Z)'$, where $N^*Z$ denotes the conormal bundle of $Z$ minus the zero section. Hence $X_{M}: \mathcal{E}'(\intm)\rightarrow \mathcal{D}'(\mcc) $ is continuous.  Here, $\mathcal{D}'(\intm), \mathcal{E}'(\intm)$ denotes the space of distributions and compactly supported distributions on $\intm$. 
 
It is known that on $\mbr^4$, the light ray transform is injective on $C_0^\infty(\mbr^4)$, see \cite{SU, Jo}, but not injective on $\mathcal{S}(\mbr^4)$ (Schwartz functions on $\mbr^4$). It is proved in \cite[Corollary 7]{Jo} that the kernel of the transform consists of  $\mathcal{S}(\mbr^4)$ functions whose Fourier transforms are supported in the time-like cone. One can obtain analogous results for $X_M$. The point is that after taking the light ray transform, time-like singularities in the functions are lost. 

To see the difference in the treatment between space-like and light-like singularities, consider the normal operator $X_{M}^*X_{M}.$ For the light ray transform on $\mbr^4$, the Schwartz kernel of the normal operator can be computed explicitly using Fourier transforms, see \cite{SU}. Let's look at the microlocal structure. The  canonical relation $C = N^*Z'$  is 
\beqq\label{eqcan}
\begin{gathered}
C = \{((y, v, \eta, w); (t, x, \tau, \xi)) \in  (T^*\mcc\backslash 0)  \times (T^*\intm \backslash 0) : y = x- t v,\ \ \eta = \xi, \\
w = t \xi|_{T_v\mbs^2},\ \ \tau = -\xi \cdot v, \ \ y\in \mbr^3, v\in \mbs^2, \eta \in \mbr^3, (t, x) \in \intm\},
\end{gathered}
\eeqq
see \cite[equation (39)]{LOSU}. In the expression of $w$, $\xi$ is regarded as a co-tangent vector to $T_v\mbs^2.$ If $\Xi = (\tau, \xi)$ is light-like, then $\xi|_{T_v\mbs^2} = 0$, see \cite[Lemma 10.1]{LOSU}. We look at the double fibration picture 
\begin{center} 
\begin{tikzpicture}
  \matrix (m) [matrix of math nodes,row sep=1.5em,column sep=1em,minimum width=1em]
  {
      & C & \\
     T^*M &  & T^*\mcc \\};
  \path[-stealth]
    (m-1-2) edge node [left] {$\pi$} (m-2-1)
    (m-1-2) edge node [right] {$\rho$} (m-2-3);
\end{tikzpicture}
\end{center}
If $\rho$ is an injective immersion, the double fibration satisfies the Bolker condition, and the normal operator $X_{M}^*\circ X_{M}$ belongs to the clean intersection calculus so that  the normal operator is a pseudo-differential operator, see  \cite{Gu0}. As shown in \cite[Lemma 10.1]{LOSU}, 
$\rho$ fails to be injective on the set $\mcl\cap C $ where
\beq
\mcl = \{(y, v, \eta, w; t, x, \Xi) \in (T^*\mcc\backslash 0) \times (T^*\intm\backslash 0): \text{$\Xi$ is light-like}\}.
\eeq
In particular, the normal operator is an elliptic pseudo-differential operator when restricted to space-like directions, see \cite{SU} and \cite{LOSU}. In general, it is proved in \cite{Wa} that the Schwartz kernel of the normal operator $X_M^*X_M$ is a paired Lagrangian distribution and a parametrix can be constructed within the framework of \cite{GU89}.  However, the picture near light-like directions is still not so clear.  We remark that Guillemin \cite{Gu} considered the structure of $X_M X_M^*$ for $2+1$ dimensional Minkowski spacetime.

\section{Solution of the Cauchy problem}\label{sec-cauchy}
We find a representation of the solution of the Cauchy problem in this section. Consider  
\beqq\label{eq-scalar0}
\begin{gathered}
 \square_c  u   = 0, \quad \text{ on }  \intm = (t_0, t_1)\times \mbr^3 \\
u = f_1, \quad \p_t u = f_2, \text{ on } \mcs_0 = \{t_0\}\times \mbr^3. 
\end{gathered}
\eeqq
The fundamental solution can be written down quite explicitly. However, it will be more convenient to look at its microlocal structure.  For \eqref{eq-scalar0}, all we need is the Fourier transform, see for example Tr\`eves \cite[Chapter VI, Section 1]{Tre}.  For general strictly hyperbolic equations, Duistermaat-H\"ormander (see \cite[Chaper 5]{Du}) constructed a parametrix for the Cauchy problem. So one can find a parametrix for \eqref{eq-scalar0} even when the equation contains lower order terms which will be used in Section \ref{sec-full}.

Let $(\tau, \xi), \xi\in \mbr^3$ be the dual variables in $T^*\intm$ to $(t, x), x\in \mbr^3$. Taking the Fourier transform of \eqref{eq-scalar0} in the $x$ variable, we get (for $t_0 = 0$)
\beq 
\begin{gathered}
\p_t^2 \hat u(t, \xi) + c^2 |\xi|^2 \hat u(t, \xi)   = 0,   \\
\hat u(0, \xi) = \hat f_1(\xi), \quad \p_t \hat u(0, \xi) = \hat f_2(\xi) . 
\end{gathered}
\eeq 
Solving this ODE, we get 
\beq
\hat u(t, \xi) = \ha e^{i t c |\xi|}(\hat f_1 + \frac{1}{ic|\xi|}\hat f_2) +\ha e^{-itc |\xi|} (\hat f_1 - \frac{1}{ic|\xi|}\hat f_2).
\eeq
Taking the inverse Fourier transform, we get 
\beqq\label{eq-cauchysol}
\begin{gathered}
u(t, x) = (2\pi)^{-3}\ha \int_{\mbr^3} e^{i(x\cdot \xi + ct |\xi|)} (\hat f_1 + \frac{1}{ic|\xi|}\hat f_2) d\xi + (2\pi)^{-3}\ha \int_{\mbr^3} e^{i(x\cdot \xi - tc |\xi|)} (\hat f_1 - \frac{1}{ic|\xi|}\hat f_2) d\xi\\
 =   (2\pi)^{-3}\int_{\mbr^3} e^{i(x\cdot \xi + ct |\xi|)}  \hat h_1(\xi) d\xi +  (2\pi)^{-3}\int_{\mbr^3} e^{i(x\cdot \xi - tc |\xi|)}  \hat h_2(\xi) d\xi \\
  = E_+ h_1 + E_- h_2, 
\end{gathered}
\eeqq
where 
\beq
\hat h_1 = \ha (\hat f_1 + \frac{1}{ic|\xi|}\hat f_2), \ \ \hat h_2 =  \ha (\hat f_1 - \frac{1}{ic|\xi|}\hat f_2).
\eeq
We see that $E_\pm$ are represented by oscillatory integrals 
\beqq\label{eq-paraE}
E_\pm(f) =(2\pi)^{-3} \int_{\mbr^3}\int_{\mbr^3} e^{i((x-y)\cdot \xi \pm ct |\xi|)}  f(y) dy d\xi.
\eeqq
The phase functions are  
$
\phi_\pm(t, x, y, \xi) = (x - y) \cdot \xi \pm ct |\xi|
$
and amplitude function $a(t, x, \xi) = 1.$ In H\"ormander's notation, we conclude that $E_\pm \in I^{-\frac 14}(\mbr^3\times\intm; (C^\pm)')$ are Fourier integral operators where the canonical relations are 
\beqq\label{eq-cano0}
C^\pm = \{(t, x, \zeta_0, \zeta'; y, \xi) \in T^*\intm \backslash 0\times T^*\mbr^3\backslash 0: y = x - c t (\pm \xi/|\xi|), \zeta'  = \xi, \zeta_0 = \pm c |\xi|\}. 
\eeqq
It suffices to regard $h_1, h_2$ as the reparametrized initial conditions for the Cauchy problem and represent $u = E_+ h_1 + E_- h_2$ in \eqref{eq-cauchysol}. Once we find $h_1, h_2$, we can easily find $f_1, f_2$ from
\beqq\label{eq-fh}
f_1 = h_1 + h_2, \ \ f_2 =  i c\lap^\ha (h_1 - h_2).
\eeqq

\section{The microlocal inversion: $c < 1$}\label{sec-local1}
For $0 < c < 1$, it is important to observe that  singularities (or the wave front set) of the solution $u$ to \eqref{eq-scalar0} are all in space-like directions for $(M, g_M).$ From the canonical relation $C^\pm$ in \eqref{eq-cano0}, we know that for $u$ in \eqref{eq-scalar0}
\beq
\WF(u) \subset  \{(t, x, \xi_0, \xi') \in T^*\intm \backslash 0: \xi_0 = \pm c |\xi'|\},
\eeq
 and $|(\xi_0, \xi')|^2_{g^*_M} = -\xi_0^2 + |\xi'|^2 = (-c^2+1)|\xi'|^2 > 0$ for $c<1$. For such $(\xi_0, \xi')$, the corresponding vector in $T\intm$ is time-like. So these singularities correspond  to trajectories of particles moving slower than photons in $(M, g_M)$. 

Now we can use the fact that in space-like directions, the normal operator $X_M^*\circ X_M$ is actually a pseudo-differential operator as shown in \cite{LOSU}.  The symbol of $\square_c$ is $p_c(\xi_0, \xi') = -\xi^2_0 + c^2 |\xi'|^2. $
Let $\chi(t)$ be a smooth  cut-off function with $\chi(t) = 1, |t|< 1$ and $\chi(t) = 0, |t|> 1/c^2$ for $c<1$. Then we define
\beq
\chi_1(\xi_0, \xi') = \chi(\frac{\xi_0^2}{c^2|\xi'|^2})
\eeq
so $\chi_1(\xi_0, \xi') = 1$ on $\{(\xi_0, \xi') \in \mbr^4: p_c(\xi_0, \xi') > 0\}$ and $\chi_1(\xi_0, \xi') = 0$ on $\Omega^{*, -}\intm$. Let $\chi_1(D)$ be the pseudo-differential operator with symbol $\chi_1$. We have

\begin{lemma}
$\chi_1(D) X_M^*\circ X_M \chi_1(D)$ is a pseudo-differential operator of order $-1$ on $\intm$. The principal symbol at $(t, x, \xi_0, \xi') \in T^*\intm$ is 
\beq
 \frac{4\pi^2}{|\xi'|}\chi_1^2(\xi_0, \xi').
 \eeq
\end{lemma}
\bpf
It follows from Theorem 2.1 of \cite{LOSU1} that $\chi_1(D) X_M^*\circ X_M \chi_1(D)$ is a pseudo-differential operator on $\intm$ with an oscillatory integral representation.  The symbol is 
\beqq\label{eq-symnormal}
\begin{gathered}
\sigma(t, x, \xi_0, \xi') = 2\pi |\mbs^1| \chi_1^2(\xi_0, \xi') |\xi'|^{-1}
\end{gathered}
\eeqq
We remark that the symbol is singular at $\xi = 0$ but this can be removed by introducing a smooth cut-off function supported near $\xi = 0$ and noticing that $|\xi|^{-1}$ is integrable near $\xi = 0$. Since it only changes $\chi_1(D) X_M^*\circ X_M \chi_1(D)$ by a smoothing operator, we will not show it for simplicity. 
\epf

Now we show that 
\begin{lemma}\label{lm-elliptic1}
The normal operator $E_+^*X_M^*\circ X_M E_+, E_-^*X_M^*\circ X_M E_-$ are  elliptic pseudo-differential operators of order $-1$ on $\mbr^3$, and $E_+^*X_M^*\circ X_M E_-$ and $E_-^*X_M^*\circ X_M E_+$ are smoothing operators on $\mbr^3.$
\end{lemma}
\bpf
First of all, we know that  $(X_M^*\circ X_M) E_+  = (\chi_1(D) X_M^*\circ X_M \chi_1(D)) E_+$ modulo a smoothing operator, thus $(X_M^*\circ X_M) E_+ \in I^{-\frac 54}(\intm\times\mbr^3; (C^+)')$ from the composition of a pseudo-differential operator and an FIO. The principal symbol is non-vanishing. We also know that $E_+^* \in I^{-\frac 14}(\intm\times\mbr^3; (C^{+, -1})')$. To compose these two operators, we would like to apply the clean composition theorem \cite[Theorem 25.2.3]{Ho4}, however, the operators are not properly supported.  But this can be justified using the oscillatory integral representation. We have (modulo a pseudo-differential operator of a lower order)
\beq
\begin{gathered}
E_+^*(X_M^*\circ X_M E_+) f(z) = (2\pi)^{-6}  \int_{\mbr^3}\int_{0}^{t_1}\int_{\mbr^3}\int_{\mbr^3}\int_{\mbr^3} e^{i((z-x)\cdot \eta - c t |\eta|)} e^{i((x-y)\cdot \xi + ct |\xi|)} a(\xi) f(y) dy d\xi dx dt d\eta\\
 = (2\pi)^{-6}  \int_{\mbr^3}\int_{0}^{t_1}\int_{\mbr^3}\int_{\mbr^3}\int_{\mbr^3} e^{i(z\cdot \eta - y\cdot \xi + x(\xi - \eta) - c t |\eta| + ct |\xi|)} a(\xi) f(y) dy d\xi dx dt d\eta\\
  = (2\pi)^{-3} \int_{\mbr^3}\int_{\mbr^3 } e^{i(z\cdot \xi - y\cdot \xi)} t_1 a(\xi) f(y) dy d\xi.
\end{gathered}
\eeq
This is a pseudo-differential operator of order $-1$ on $\mbr^3.$ The same proof works for the minus sign.

 To see that  $E_+^*X_M^*\circ X_M E_-$ is smoothing, we just need to observe that the canonical relations $C^+, C^-$ in \eqref{eq-cano0} are disjoint. So a wave front analysis using e.g.\ \cite[Theorem 1.3.7]{Du} tells that the operator is smoothing.  
\epf
We finished the proof but we mention the following alternative argument. Essentially, we want to consider the operator $E_+$ for fixed $t$, denoted by $E_+(t)$. We know that $E_+(t): \mathcal{E}'(\mbr^3)\rightarrow \mathcal{D}'(\mbr^3)$ 
is a Fourier integral operator 
\beq
E_+(t)f(x) =(2\pi)^{-3} \int_{\mbr^3}\int_{\mbr^3 } e^{i((x-y)\cdot \xi + ct |\xi|)}  f(y) dy d\xi 
\eeq
with canonical relation 
$
C_t = \{(y, \eta; x, \xi)\in T^*\mbr^3\backslash 0 \times T^*\mbr^3\backslash 0 : y = x + ct \xi/|\xi|, \xi = \eta\}. 
$
Then $E_+(t)\in I^{0}(\mbr^3\times \mbr^3; C_t')$ is properly supported. The canonical relation $C_t$ is a graph of a symplectic transformation, thus the composition $E_+^{*}(t)E_+(t)$ is a pseudo-differential operator of order $0$ on $\mbr^3$. In our case, $E_+^{*}(t)X_M^*X_M E_+(t)$ is a pseudo-differential operator of order $-1$ and the symbols are smooth in $t \in [t_0, t_1]$. Finally, integrating the symbols in $t$ produces a symbol and we get the result. 

Now we construct a parametrix for the transform.  
\begin{prop}\label{prop-para1}
For $c< 1,$ there exist  operators $A_1, A_2$ such that 
\beq
A_1 X_M f = f_1 + R_1 f_1 + R_1'f_2, \ \ A_2 X_M f = f_2 + R_2 f_1 + R_2'f_2
\eeq
where $R_1, R_2, R_1', R_2'$ are smoothing operators and $A_i = \tilde A_i \circ X_M^*, i = 1, 2$ in which $\tilde A_i$ are Fourier integral operators.  
\end{prop}
\bpf
First, we represent $f = E_+ h_1 + E_- h_2$ and write 
\beqq\label{eq-xmf}
X_M f = X_M E_+ h_1 + X_M E_- h_2.
\eeqq
We apply $E_+^* X_M^*$ to get 
\beq
E_+^* X_M^*X_M f =E_+^* X_M^* X_M E_+ h_1 + E_+^* X_M^*X_M E_- h_2 = E_+^* X_M^* X_M E_+ h_1 + R_1 h_2.
\eeq 
Since $E_+^* X_M^* X_M E_+$ is an elliptic pseudo-differential operator of order $-1$, we can find a parametrix $B_+$ which is a pseudo-differential operator of order $1$ on $\mbr^3$ and  
\beq
B_+\circ E_+^* X_M^*X_M f = h_1 +  R_1 h_1 + R_1' h_2
\eeq
where $R_1, R_1'$ are smoothing. We repeat the argument for the minus sign. Apply $E_-^*X_M^*$ to \eqref{eq-xmf}, we get 
\beq
E_-^* X_M^*X_M f =E_-^* X_M^* X_M E_+ h_1 + E_-^* X_M^*X_M E_- h_2 = E_-^* X_M^* X_M E_- h_2 + R_2 h_2.
\eeq 
Apply the parametrix $B_-$ for $E_-^* X_M^* X_M E_-$ and we get 
\beq
B_-\circ E_-^* X_M^*X_M f = h_2+  R_2 h_1 + R_2' h_2.
\eeq
Finally, we get 
\beq
\begin{gathered}
f_1 + R_1 f_1 + R_2 f_2 = (B_+\circ E_+^* + B_- \circ E_-^*) X_M^*X_M f    \\
\text{ and }  
f_2 + R_1' f_1 + R_2' f_2 = ic\lap^\ha (B_+\circ E_+^* + B_- \circ E_-^*) X_M^*X_M f 
\end{gathered}
\eeq
as claimed. We set $\tilde A_1 = B_+\circ E_+^* + B_- \circ E_-^* $ which is a sum of two FIOs in $I^{3/4}(\intm\times\mbr^3; (C^{+, -1})')$ and $I^{3/4}(\intm\times\mbr^3; (C^{-, -1})')$, and $\tilde A_2 = ic\lap^\ha (B_+\circ E_+^* + B_- \circ E_-^*)$ which is a sum of two FIOs in $I^{7/4}(\intm\times\mbr^3; (C^{+, -1})')$ and $I^{7/4}(\intm\times\mbr^3; (C^{-, -1})')$. This completes the proof. 
\epf

For convenience, we formulate a microlocal inversion result for determining $f$.
\begin{cor}
For $c < 1$, there exist  operators $A$ such that 
\beq
A  X_M f = f  + R_1 f_1 + R_2f_2, 
\eeq
where $R_1, R_2$ are smoothing operators. 
\end{cor}
\bpf
Again, we simply solve the wave equation \eqref{eq-scalar0} using the parametrix. In fact, it is easier to use $h_1, h_2$. 
\beq
\begin{gathered}
f = E_+ h_1 + E_-h_2   
= E_+ B_+ \circ E_+^* X_M^*X_M f  + E_-B_-\circ E_-^* X_M^*X_M f  + \tilde R_1 h_1 +\tilde R_2 h_2 \\
 = (E_+ B_+ \circ E_+^* + E_-B_-\circ E_-^*) X_M^*X_M f + R_1 f_1 + R_2 f_2
\end{gathered}
\eeq
as claimed, where $\tilde R_1, \tilde R_2, R_1, R_2$ are smoothing operators and $A = (E_+ B_+ \circ E_+^* + E_-B_-\circ E_-^*) X_M^*.$
\epf

\section{The microlocal inversion: $c = 1$}\label{sec-local2}
For $c = 1$, the singularities of the solutions of \eqref{eq-scalar0} are all in light-like directions. As explained in the end of Section \ref{sec-light}, the Schwartz kernel of $X_M^* \circ X_M$ is more complicated and the previous argument does not work directly. We will take a different approach by considering the composition $X_M\circ E_\pm$. Let $\varphi$ be a smooth function on $\mbs^2$, and $I^\varphi$ be the integration operator on $C^\infty(\mbr^3\times \mbs^2)$ defined by
\beq
I^\varphi f(y) = \int_{\mbs^2} \varphi(v) f(y, v)dv.
\eeq  
Then we consider the composition $I^\varphi\circ X_M \circ E_\pm$ as an operator from $C^\infty(\mcs_0)$ to $C^\infty(\mcs_0)$.  For technical reasons, we introduce a smooth cut-off function. For $\eps>0$ small, let $\chi_\eps(t)$ be a smooth cut-off function on $\mbr$ such that $\chi_\eps(t) = 1$ for $2\eps < t < t_1 - 2\eps$ and $\chi_\eps(t) = 0$ for $t<\eps$ and $t > t_1 - \eps.$ We prove
\begin{prop}\label{prop-compose}
$K_\pm \doteq I^\varphi X_M \chi_\eps E_\pm \in \Psi^{-1}(\mcs_0)$ are pseudo-differential operators of order $-1$ with complete symbol $k_\pm(\xi), \xi\in \mbr^3\backslash 0$ and the principal symbols are given by 
\beq
\begin{gathered}
k_{+, -1}(\xi) = 2\pi i  c_\eps |\xi|^{-1}  \varphi(-\xi/|\xi|), \quad 
k_{-, -1}(\xi) = -2\pi i  c_\eps |\xi|^{-1}\varphi(\xi/|\xi|), \\
\text{where } c_\eps = \int_0^{t_1}t^{-1} \chi_\eps(t)dt 
\end{gathered}
\eeq
\end{prop} 
\bpf
We start with $K_+.$ We recall from  \eqref{eqlit} that 
\beq
\begin{gathered}
X_{M}f(y, v)   = (2\pi)^{-3} \int_{\mbr^3}\int_{\mbr^3}\int_{0}^{t_1} e^{i ( (y-x)\cdot \eta + t v\cdot \eta)} f(t, x)  dt dx d\eta
\end{gathered}
\eeq 
and from Section \ref{sec-cauchy} that
\beq
E_+(f)(t, x) =(2\pi)^{-3} \int_{\mbr^3}\int_{\mbr^3} e^{i((x-z)\cdot \xi +  t |\xi|)}  f(z) dz d\xi. 
\eeq
Consider  the oscillatory integral integral representation of the Schwartz kernel $K_+$  
\beqq\label{eq-XME}
 \begin{gathered}
K_+(y, z)   = (2\pi)^{-6} \int_{\mbs^2} \int_{\mbr^3}\int_{\mbr^3}\int_{0}^{t_1} \int_{\mbr^3 } e^{i ( (y - x)\cdot \eta + t v\cdot \eta+ (x-z)\cdot \xi +   t |\xi|)} \varphi(v) \chi_\eps(t) d\xi dt dx d\eta dv   
\end{gathered}
\eeqq
In this case, the oscillatory integral can be computed explicitly. But before we proceed with the calculation, we examine the phase function 
\beq
\phi(y, z, \xi, t; \eta, x, v) = (y - x)\cdot \eta + t v\cdot \eta+ (x-z)\cdot \xi +   t |\xi|
\eeq
Consider $\phi$ in $\eta, x, v$ variables. We have
\beq
\phi_\eta = y - x + tv, \quad \phi_x = \xi - \eta,  \quad\phi_v = t\eta|_{T_v\mbs^2},
\eeq
so  the critical points are given by  
\beq
\xi = \eta, \quad v = \pm \xi/|\xi|, \quad x = y - t \xi/|\xi|
\eeq
Here, we remark that $t\xi|_{T_v\mbs^2} = 0$ implies that $\xi$ is parallel to $v$ so $v = \pm \xi/|\xi|$.  Also, we have 
\beq
 \p^2_{(\eta, x, v)} \phi  = \begin{pmatrix}
0 & -\id &t \\
  \id &  0&0\\
 \ast & 0 & \ast 
\end{pmatrix}
\eeq
 To compute $\ast$, we introduce local coordinates on $\mbs^2$ near the critical point. By using an orthogonal transformation, we can assume that $\xi/|\xi| = (0, 0, 1).$ We use $v = (v_1, v_2, \pm \sqrt{1 - v_1^2 -v_2^2})$ near $\pm \xi/|\xi|$ where $v_1^2 + v_2^2 <1$ . Then we have
\beq
\p_v \phi = \p_{(v_1, v_2)}(t v\cdot \eta)=  t
\begin{pmatrix}
\eta_1 \pm   \eta_3 \frac{-v_1}{\sqrt{1 - v_1^2 - v_2^2}}\\
\eta_2 \pm  \eta_3 \frac{-v_2}{\sqrt{1 - v_1^2 - v_2^2}}
\end{pmatrix}.
\eeq
On the set of critical point, $v = \pm (0, 0, 1)$ and $\eta = (0, 0,  |\xi|)$. We observe that $\p_v \phi = 0.$ Next,  
\beq
\p_{\eta}(\frac{\p\phi}{\p v}) = t\begin{pmatrix}
1          & 0 &  \pm \frac{-v_1}{\sqrt{1 - v_1^2 - v_2^2}}  \\
0 &  1  & \pm \frac{-v_2}{\sqrt{1 - v_1^2 - v_2^2}}
\end{pmatrix}
 \text{ and }
\p_{v}(\frac{\p \phi}{\p v}) =  \pm t \eta_3 \begin{pmatrix}
\frac{-1 + v_2^2}{(1 - v_1^2 - v_2^2)^{\frac 32}} & \frac{-v_1 v_2}{(1 - v_1^2 - v_2^2)^{\frac 32}} \\
 \frac{-v_1 v_2}{(1 - v_1^2 - v_2^2)^{\frac 32}}   &  \frac{-1 + v_1^2}{(1 - v_1^2 - v_2^2)^{\frac 32}}
\end{pmatrix}.  
\eeq
On  critical points,  
 \beq
\begin{gathered}
\p_{v}(\frac{\p \phi}{\p v})  
 = \pm t|\xi| \begin{pmatrix}
-1&0\\
0 &  -1
\end{pmatrix}.
 \end{gathered}
\eeq
This shows that the phase function is non-degenerate in $\eta, x, v$. We can apply stationary phase argument so the phase becomes 
\beq
\begin{gathered}
\phi(y, z, \xi, t) = (y - z)\cdot \xi + 2t |\xi|  \text{ when } v = \xi/|\xi|\\
\phi(y, z, \xi, t) = (y - z)\cdot \xi   \text{ when } v = -\xi/|\xi|
\end{gathered}
\eeq
Finally, after integrating in $t$, we will get a pseudo-differential operator. This will be shown explicitly in the follows. 

First, in \eqref{eq-XME}, we integrate in $x, \eta$  to get 
\beq 
 \begin{gathered}
K_+(y, z)  
 = (2\pi)^{-3} \int_{\mbs^2}  \int_{0}^{t_1}\int_{\mbr^3} e^{i ( y\cdot \xi + t v\cdot \xi -z \cdot \xi + t |\xi|)} \varphi(v)\chi_\eps(t) d\xi dt   dv 
\end{gathered}
\eeq 
Consider the integral in $v$. For $t$ non-zero, the $v$ integral is non-degenerate with stationary points at $v=\pm \xi/|\xi|$. Applying stationary phase argument see e.g. \cite[Lemma 1.2]{Mel}, we get 
\beqq\label{eq-XME2}
 \begin{gathered}
K_+(y, z) 
= (2\pi)^{-3}  \int_{0}^{t_1}\int_{\mbr^3}  e^{i ( (y-z) \cdot \xi + 2 t |\xi|)} (\varphi(\xi/|\xi|) + \varphi^+(t, \xi))\chi_\eps(t) (t|\xi|)^{-1} e^{-\frac{1}{2}i\pi}(2\pi)   d\xi dt \\
+ (2\pi)^{-3}  \int_{0}^{t_1}\int_{\mbr^3}  e^{i (y-z) \cdot \xi} (\varphi(-\xi/|\xi|) + \varphi^-(t, \xi)) \chi_\eps(t) (t|\xi|)^{-1} e^{\frac{1}{2}i\pi}(2\pi) d\xi dt \\
 =  (2\pi)^{-3} \int_{\mbr^3}  e^{i (y-z) \cdot \xi  } k_+(\xi) d\xi 
\end{gathered}
\eeqq
where $\varphi^\pm$ come from the stationary phase argument and they have asymptotic expansions
\beqq\label{eq-asympphi}
\varphi^\pm(t, \xi) \sim \sum_{k = 1}^\infty a_k^\pm(\xi/|\xi|) (t|\xi|)^{-k}
\eeqq
in which $a_k^\pm$ are smooth functions on $\mbs^2$. Also, 
\beqq\label{eq-kplus}
\begin{gathered}
k_+(\xi) =  
 + 2\pi i|\xi|^{-1} \varphi(-\xi/|\xi|) \int_{0}^{t_1} t^{-1}(1 + \varphi^-(t, \xi))\chi_\eps(t)dt\\
 - 2\pi i|\xi|^{-1}  \varphi(\xi/|\xi|) \int_0^{t_1}e^{2it|\xi|}t^{-1}(1 +\varphi^+(t, \xi))\chi_\eps(t)dt  
\end{gathered}
\eeqq
The second integral in $t$ is $O(|\xi|^{-\infty})$ for $|\xi|$ large because $t$ is away from $0$ and $\chi_\eps$ is smooth. For the first integral, the integral of each asymptotic term of $\varphi^-$ in \eqref{eq-asympphi} in $t$ is finite. Thus $k_+(\xi)$ is a symbol of order $-1$ and the leading order term is 
\beq
k_{+, -1}(\xi) = 2\pi i |\xi|^{-1}  \varphi(-\xi/|\xi|) \int_0^{t_1}t^{-1} \chi_\eps(t)dt.  
\eeq 
This shows that $K_+$ in \eqref{eq-XME2} is a pseudo-differential operator of order $-1$ on $\mbr^3.$ 

For $K_-$, the calculation is similar and we look for the symbol. 
\beq 
 \begin{gathered}
K_-(y, z)  = (2\pi)^{-3}  \int_{t_0}^{t_1}\int_{\mbr^3} e^{i ( y\cdot \xi + t v\cdot \xi -z \cdot \xi - t |\xi|)} \chi_\eps(t) d\xi dt   dv\\
 = - i(2\pi)^{-2}  \int_{t_0}^{t_1}\int_{\mbr^3}  e^{i  (y-z) \cdot \xi} (t|\xi|)^{-1}  (\varphi(\xi/|\xi|) + \tilde\varphi^+(t, \xi)) \chi_\eps(t)d\xi dt \\
+i (2\pi)^{-2}  \int_{t_0}^{t_1}\int_{\mbr^3}  e^{i ((y-z) \cdot \xi - 2t|\xi|)} (t|\xi|)^{-1}(\varphi(-\xi/|\xi|) + \tilde \varphi^-(t, \xi))\chi_\eps(t) d\xi dt  \\
 =   (2\pi)^{-3}\int_{\mbr^3}  e^{i (y-z) \cdot \xi  } k_-(\xi) d\xi 
\end{gathered}
\eeq 
where $\tilde\varphi^\pm$ have similar asymptotic expansion as \eqref{eq-asympphi}, and $k_-(\xi)$ is given by 
\beqq\label{eq-kminus}
\begin{gathered}
k_-(\xi) = - 2\pi i |\xi|^{-1} \int_0^{t_1}(\varphi(\xi/|\xi|) + \tilde\varphi^+(t, \xi))  t^{-1}   \chi_\eps(t)dt  \\
 +  2\pi i |\xi|^{-1} \int_{0}^{t_1} e^{-2it|\xi|}(\varphi(-\xi/|\xi|) + \tilde\varphi^-(t, \xi))  t^{-1}\chi_\eps(t)dt 
\end{gathered}
\eeqq
This  is a symbol of order $-1$ and the leading order term is 
\beq
k_{-, -1}(\xi) = -2\pi i |\xi|^{-1}\varphi(\xi/|\xi|)  \int_0^{t_1}   t^{-1}   \chi_\eps(t)dt 
\eeq
This completes the proof of the proposition.  
\epf

Next we discuss what needs to be changed when the smooth cut-off function $\chi_\eps$ is replaced by the characteristic function $\chi_{[\eps, t_1]}$ of the interval $[\eps, t_1]$ in $\mbr$.    
All the calculations in Proposition \ref{prop-compose} hold up to \eqref{eq-kplus} which is now  
\beqq\label{eq-kplus1}
\begin{gathered}
k_+(\xi) =  
 + 2\pi i|\xi|^{-1} \varphi(-\xi/|\xi|) \int_{\eps}^{t_1} t^{-1}(1 + \varphi^-(t, \xi)) dt\\
 - 2\pi i|\xi|^{-1}  \varphi(\xi/|\xi|) \int_{\eps}^{t_1}e^{2it|\xi|}t^{-1}(1 +\varphi^+(t, \xi)) dt  
\end{gathered}
\eeqq
The first integral, denoted by $I_1$ below, still gives a symbol of order $-1$. For the second integral denoted by $I_2$ below, integration by parts gives
\beq 
\begin{gathered}
I_2(\xi) = - 2\pi i|\xi|^{-1}  \varphi(\xi/|\xi|) \{ \frac{1}{2i|\xi|} (e^{2it_1|\xi|} t_1^{-1} (1 + \varphi^+(t_1, \xi))) -  \frac{1}{2i|\xi|} e^{2i \eps |\xi|} \eps^{-1} (1 + \varphi^+(\eps, \xi)) \\
-   \frac{1}{2i|\xi|} \int_{\eps}^{t_1}e^{2it|\xi|}\frac{d}{dt}[t^{-1}(1 +\varphi^+(t, \xi)) ]dt \}
\end{gathered}
\eeq 
We can repeat the integration by parts and get 
\beq
I_2(\xi) =  e^{2it_1|\xi|} a(\xi) + e^{2i\eps |\xi|}b(\xi)
\eeq
where $a(\xi), b(\xi)$ are symbols of order $-2.$ 
Using these in \eqref{eq-XME2}, we get 
\beq 
 \begin{gathered}
K_+(y, z)  
 =  (2\pi)^{-3} \int_{\mbr^3}  e^{i (y-z) \cdot \xi  } I_1(\xi) d\xi +  (2\pi)^{-3} \int_{\mbr^3}  e^{i (y-z) \cdot \xi  + 2it_1|\xi|} a(\xi) d\xi \\
 +  (2\pi)^{-3} \int_{\mbr^3}  e^{i (y-z) \cdot \xi  + 2i\eps |\xi|} b(\xi) d\xi.
\end{gathered}
\eeq 
Thus, we can write $K_+ = K_+^0 + K_+^{\eps} + K_+^{t_1}$ where $K_+^0 \in \Psi^{-1}(\mbr^3)$, and $K_+^\eps \in I^{-2}(\mbr^3, \mbr^3; C_\eps), K_+^{t_1}\in I^{-2}(\mbr^3, \mbr^3; C_{t_1})$ are Fourier integral operators of order $-2.$ The canonical relation $C_\eps, C_{t_1}$ can be described as follows. For $\alpha \in \mbr$, we define 
\beq
\begin{gathered}
C_\alpha = \{(y, \eta, z, \zeta) \in T^*\mbr^3\backslash 0\times T^*\mbr^3\backslash 0: y = z + 2 \alpha \xi/|\xi|,  \xi = \eta \}.
\end{gathered}
\eeq
We see that $C_\alpha$ is a graph of a canonical transformation, see \cite[Section 25.3]{Ho4}. The same argument shows that $K_-$ is also a sum of $K_-^0 \in  \Psi^{-1}(\mbr^3)$ and $K_-^\eps \in I^{-2}(\mbr^3, \mbr^3; C_{-\eps}), K_-^{t_1} \in I^{-2}(\mbr^3, \mbr^3; C_{-t_1})$. 
\\

Now we are ready to obtain a parallel result  of Proposition \ref{prop-para1} about the microlocal inversion. 
\begin{prop}\label{prop-para2}
For $c = 1$ and any $N\in \mbn$, there exist  operators $A_1, A_2$ such that 
\beq
A_1 X_M \chi_{[\eps, t_1]} f = h_1 + R_1h_1 + R_1'h_2, \ \ A_2 X_M \chi_{[\eps, t_1]} f = h_2 + R_2 h_1 + R_2'h_2
\eeq
where $h_1, h_2$ are defined in Section \ref{sec-cauchy} and $R_1, R_1', R_2, R_2' \in I^{-N}(\mbr^3, \mbr^3; C^N_{\eps, t_1})$ which  is  the $N$-fold composition of elements in $I^{-1}(\mbr^3, \mbr^3; C_{\pm \eps})$ and  $I^{-1}(\mbr^3, \mbr^3; C_{\pm t_1}) $, more explicitly
\beq
\begin{gathered}
I^{-N}(\mbr^3, \mbr^3; C^N_{\eps, t_1}) = \{A_1 \circ A_2\cdots A_N: A_i \in I^{-1}(\mbr^3, \mbr^3; C_{\eps}) + I^{-1}(\mbr^3, \mbr^3; C_{t_1})\\
+ I^{-1}(\mbr^3, \mbr^3; C_{-\eps}) + I^{-1}(\mbr^3, \mbr^3; C_{-t_1}) \}.
\end{gathered}
\eeq
\end{prop}
\bpf
We divide the proof in two steps. 

{\bf Step 1:} Let's replace $\chi_{[\eps, t_1]}$ with the smooth cut-off $\chi_\eps$ as in Proposition \ref{prop-compose} and see how to get $h_1, h_2$ using Proposition \ref{prop-compose}. We  write 
\beq 
X_M \chi_\eps f = X_M\chi_\eps E_+ h_1 + X_M \chi_\eps E_- h_2.
\eeq 
Let $\varphi$ be a smooth function on $\mbs^2$. Applying $I^\varphi$ we get 
\beqq\label{eq-qp1}
I^\varphi X_M\chi_\eps f = I^\varphi  X_M \chi_\eps E_+ h_1 +  I^\varphi  X_M \chi_\eps E_- h_2  = K^{\varphi, +} h_1 + K^{\varphi, -}h_2
\eeqq   
where we added $\varphi$ to the notation of $K_\pm$ to emphasize the dependency because we will choose different $\varphi$ below. 

First, let $\varphi_1 = 1$. From Proposition \ref{prop-compose}, we see that $K^{\varphi_1}_\pm \in \Psi^{-1}(\mbr^3)$ and the principal symbols  are given by 
\beq
\begin{gathered}
k^{\varphi_1}_{+, -1}(\xi) 
= -k^{\varphi_1}_{-, -1}(\xi) =  2\pi i c_\eps |\xi|^{-1}. 
\end{gathered}
\eeq  
We let $Q^1_+$ be a parametrix of $K^{\varphi_1}_+$  and get 
\beqq\label{eq-q1}
\begin{gathered}
Q^{1}_+ I^{\varphi_1} X_M\chi_\eps f  =  h_1 + Q^1_+ K^{\varphi_1}_-h_2 + R_1 h_1 
\end{gathered}
\eeqq
where $R_1, R_2$ are smoothing operators.  From the composition of pseudo-differential operators, we know that  
$Q^1_+ K^{\varphi_1}_- \in \Psi^{0}(\mbr^3)$ with principal symbol equal to $-1.$

Next, we change the function $\varphi$. Ideally, we will take an odd function $\varphi(-v) = -\varphi(v)$ but then $\varphi$ vanishes somewhere on $\mbs^2$ so we proceed as follows. Let $x = (x_1, x_2, x_3)$ be the coordinate for $\mbr^3$.  For $\delta >0$, let $\mcu_k = \{v: v = (x_1, x_2, x_3), \|x\| = 1, |x_k| > \delta/2\}, k = 1, 2, 3.$ For $\delta$ sufficiently small, $\mcu_k, k = 1, 2, 3$ form an open covering of $\mbs^2.$ Let $\chi_k(v), k = 1, 2, 3$ be a partition of unity subordinated to this covering and $\chi_k(v) = 1$ on $\mcv_k = \{v: v = (x_1, x_2, x_3), \|x\| = 1, |x_k| > \delta\}, k = 1, 2, 3.$ Here, by possibly taking $\delta$ smaller, we can assume that $\mcv_k$ also form an open covering of $\mbs^2.$ 
For $v\in \mbs^2$, we let 
\beq
\varphi_{2, k}(v) =  \chi_k(x) x_k + 2, \quad k = 1, 2, 3
\eeq
Then $\varphi_{2}(v)\neq 0$ and $\varphi_{2, k}(-v) - \varphi_{2}(v, k) \neq 0$ for $v\in \mcu_k.$   
From Proposition \ref{prop-compose}, we know that $K^{\varphi_{2, k}}_\pm \in \Psi^{-1}(\mbr^3)$ with principal symbols
\beq
\begin{gathered}
k^{\varphi_{2, k}}_{+, -1}(\xi)   =   2\pi i  c_\eps|\xi|^{-1} \varphi_{2, k}(-\xi/|\xi|), \quad 
k^{\varphi_{2, k}}_{-, -1}(\xi) =  -2\pi i c_\eps |\xi|^{-1} \varphi_{2, k}(\xi/|\xi|).
\end{gathered}
\eeq  
We consider $k = 1$ in the follows as the other cases are similar. 
Let  $Q^{2, 1}_{+}$ be a parametrix for $K^{\varphi_{2, 1}}_+$. We get 
\beq
\begin{gathered}
Q^{2, 1}_+ I^{\varphi_{2, 1}} X_M\chi_\eps f  =  h_1 + Q^{2, 1}_+ K^{\varphi_{2, 1}}_-h_2 + R_3 h_1
\end{gathered}
\eeq
where $R_3$ is a smoothing operator, and $Q^{2, 1}_+ K^{\varphi_{2, 1}}_-\in \Psi^{0}(\mbr^3)$ with principal symbol 
\beqq\label{eq-symcomp}
\sigma_0(Q^{2, 1}_+ K^{\varphi_{2, 1}}_-)(x, \xi)  = -\frac{\varphi_{2, 1}(\xi/|\xi|)}{\varphi_{2, 1}(-\xi/|\xi|)} \neq -1 
\eeqq 
when $\xi/|\xi|\in \mcu_1.$
Now we consider 
\beq
\begin{gathered}
Q^{1}_+ I^{\varphi_1} X_M\chi_\eps f  - Q^{2, 1}_+ I^{\varphi_{2, 1}} X_M\chi_\eps f
=  ( Q^1_+ K^{\varphi_1}_-  - Q^{2, 1}_+ K^{\varphi_{2, 1}}_-)h_2 + R_{1} h_1+  R_{2} h_1 - R_3 h_1
\end{gathered}
\eeq 
We observe that $A = Q^1_+ K^{\varphi_1}_-  - Q^{2, 1}_+ K^{\varphi_{2, 1}}_-$ is a pseudo-differential operator of order $0$ and the principal symbol does not vanish on $\mcu_1.$ Let $\tilde\chi_k, k = 1, 2, 3$ be a smooth partition of unity subordinated to $\mcv_k$. Then $\chi_1\tilde \chi_1 = \tilde \chi_1$. Let $B_1$ be a pseudo-differential operator of order $0$ with principal symbol $\sigma_0(B_1)(\xi) = \tilde \chi_1(\xi/|\xi|)/\sigma_0(A)(\xi)$. We can improve $B_1$ to a parametrix for $A$ so that $B_1\circ A = \tilde\chi_1(D) + R_4$ with $R_4$ smoothing. So we get 
\beq
\begin{gathered}
 B_1(Q^{1}_+ I^{\varphi_1} X_M\chi_\eps    - Q^{2, 1}_+ I^{\varphi_{2, 1}} X_M\chi_\eps )f = \tilde \chi_1(D)h_2 + R_3 h_2 + R_4 h_1
\end{gathered}
\eeq
where by abusing notations, $R_3, R_4$ are smoothing operators. We can repeat the construction for $k = 2, 3$ to get the corresponding $B_2, B_3 \in \Psi^0(\mbr^3)$. Then we arrive at 
\beqq\label{eq-q2}
\begin{gathered}
 \sum_{k = 1}^3B_k(Q^{1}_+ I^{\varphi_1} X_M\chi_\eps    - Q^{2, k}_+ I^{\varphi_{2, k}} X_M\chi_\eps )f =  h_2 + R_5 h_2 + R_6 h_1
\end{gathered}
\eeqq
with $R_5, R_6$ smoothing. This gives $A_2 = \sum_{k = 1}^3B_k(Q^{1}_+ I^{\varphi_1}     - Q^{2, k}_+ I^{\varphi_{2, k}} )$ so that $A_2 X_M \chi_{\eps} f = h_2 + R_6 h_1 + R_5h_2$. For $A_1$, we can use \eqref{eq-q1} and \eqref{eq-q2} to get 
\beq
Q^{1}_+ I^{\varphi_1} X_M\chi_\eps f  =  h_1 + Q^1_+ K^{\varphi_1, -} A_2 X_M\chi_\eps f + R_5' h_1 + R_6' h_2
\eeq
where $R_5', R_6'$ are smoothing operators. So we obtain $A_1 = Q^{1}_+ I^{\varphi_1} - Q^1_+ K^{\varphi_1}_- A_2$ so that $A_1 X_M \chi_{\eps} f = h_1 + R_5'h_1 + R_6'h_2$.

{\bf Step 2:} Now we deal with the characteristic function $\chi_{[\eps, t_1]}$. We start with 
\beq 
X_M \chi_{[\eps, t_1]} f = X_M \chi_{[\eps, t_1]} E_+ h_1 + X_M  \chi_{[\eps, t_1]} E_- h_2.
\eeq  
Applying $I^\varphi$, we get 
\beq 
I^\varphi X_M \chi_{[\eps, t_1]}  f =  K^{\varphi}_+ h_1 + K^{\varphi}_- h_2
\eeq   
where $K^{\varphi}_\pm  = I^\varphi X_M \chi_{[\eps, t_1]} E_\pm$. 
According to the arguments after Proposition \ref{prop-compose}, we can write the above as 
\beqq\label{eq-invert2} 
\begin{gathered}
I^\varphi X_M \chi_{[\eps, t_1]}  f =  (K^{\varphi, 0}_+ + K^{\varphi, \eps}_+ + K^{\varphi, t_1}_+ )  h_1 
 + (K^{\varphi, 0}_- + K^{\varphi, \eps}_- + K^{\varphi, t_1}_- )h_2\\
 \end{gathered}
\eeqq
where $K^{\varphi, 0}_\pm \in \Psi^{-1}(\mbr^3)$, $K^{\varphi, \eps}_\pm \in I^{-2}(\mbr^3, \mbr^3; C_{\pm \eps})$ and $K^{\varphi,  t_1}_\pm \in I^{-2}(\mbr^3, \mbr^3; C_{\pm t_1})$. As in Step 1, we can apply pseudo-differential operators $Q_+^1, Q_+^{2,k}, k = 1, 2, 3$ to \eqref{eq-invert2}. The arguments for $K^{\varphi, 0}_\pm$ are the same as before. As for $K^{\varphi, \eps}_\pm, K^{\varphi, t_1}_\pm$, we notice that the composition $Q_+^1 K^{\varphi, j}_\pm, Q_+^{2,k}K^{\varphi, j}_\pm, k = 1, 2, 3, j = \eps, t_1$ are all Fourier integral operators of order $-1$ with canonical relation $C_{\pm \eps}$ or $C_{\pm t_1}$. Therefore, using the same $A_1, A_2$ in Step 1, we obtain 
\beqq\label{eq-am1}
A_1 X_M \chi_{[\eps, t_1]} f = h_1 + R_1h_1 + R_1'h_2, \ \ A_2 X_M \chi_{[\eps, t_1]} f = h_2 + R_2 h_1 + R_2'h_2
\eeqq
where $R_1, R_1', R_2, R_2' \in I^{-1}(\mbr^3, \mbr^3; C_\eps) + I^{-1}(\mbr^3, \mbr^3; C_{t_1}) + I^{-1}(\mbr^3, \mbr^3; C_{-\eps}) + I^{-1}(\mbr^3, \mbr^3; C_{-t_1})$. 

Finally, we improve the remainder term using the Neumann series. 
We write \eqref{eq-am1} in matrix form 
\beq
\begin{pmatrix}
A_1 X_M \chi_{[\eps, t_1]} f \\
A_2 X_M \chi_{[\eps, t_1]} f
\end{pmatrix}
 = \id 
 \begin{pmatrix}
h_1\\
h_2
\end{pmatrix}
+R
\begin{pmatrix}
h_1\\
h_2
\end{pmatrix}, \ \ R = \begin{pmatrix}
R_1 & R_1'\\
R_2 & R_2'
\end{pmatrix}.
\eeq
For $N \in \mbn$, we let $W = \sum_{n = 0}^{N-1} (-R)^n$ and  get 
\beq
W \begin{pmatrix}
A_1 X_M \chi_{[\eps, t_1]} f \\
A_2 X_M \chi_{[\eps, t_1]} f
\end{pmatrix}
 = \id 
 \begin{pmatrix}
h_1\\
h_2
\end{pmatrix}
+R^N
\begin{pmatrix}
h_1\\
h_2
\end{pmatrix} 
\eeq
Because $R_1, R_1', R_2, R_2'$ are FIOs of the canonical graph type, we can apply the composition result in \cite[Section 25.3]{Ho4} to conclude that the terms in $R^N$ belongs to  $I^{-N}(\mbr^3, \mbr^3; C^N_{\eps, t_1})$. Finally, we set 
\beq
\begin{pmatrix}
\tilde A_1\\
\tilde A_2
\end{pmatrix} = W
\begin{pmatrix}
A_1  \\
A_2 
\end{pmatrix}
\eeq 
Changing notations of $\tilde A_1, \tilde A_2$ to $A_1, A_2$ finishes the proof.  
\epf

\section{The stable determination}\label{sec-full}
We prove Theorem \ref{thm-main1}, starting with the injectivity of the light ray transform. It is known, see for instance \cite{SU, Jo}, that the light ray transform on $\mbr^{n+1}$ is injective on $C_0^\infty$ functions. This also holds for $L^1_{\comp}$ functions and the proof is similar, see \cite{SU}. 
\begin{theorem}\label{thm-inj}
Suppose $f\in L^1_{\comp}(\mbr^{n+1}), n \geq 2$ and $X_{\mbr^{n+1}}f = 0$.   
Then $f = 0.$ 
\end{theorem}
\bpf 
For $f\in L^1_{\comp}(\mbr^{n+1})$, the Fourier transform $\hat f$ is analytic. Let $\theta\in \mbs^{n-1}$ and $\Theta = (1, \theta)$ be a light-like vector. Let $z = (s, y + s\theta) \in \mbr^{n+1}, s\in \mbr, y\in \mbr^n$. We   parametrize the light ray transform as 
\beq
X_{\mbr^{n+1}} f(z, \Theta) = \int_{\mbr} f(t, y + t\theta)d t. 
\eeq
From the standard Fourier Slice Theorem for geodesic ray transforms on $\mbr^{n+1}$, we get 
\beq
\hat f(\zeta) = \int_{\Theta^\perp} e^{-iy\cdot \zeta} X_{\mbr^{n+1} } f(z,  \Theta)  dS_z
\eeq
where the integration is over  the hyperplane $\Theta^\perp$ perpendicular to $\Theta$ with respect to the Euclidean inner product in $\mbr^{n+1}$ and   $\zeta = (\tau, \xi)\in \mbr^{n+1}, \xi\in \mbr^n, \xi\neq 0$ is perpendicular to $\Theta$. We notice that if  $|\tau|\leq |\xi|$, then there is a null vector $(1, \theta)$ which is Euclidean orthogonal to $\zeta$. Actually, $\tau + \theta\cdot \xi = 0$ so $\theta\cdot (\xi/|\xi|) = -\tau/|\xi| \in [-1, 1]$ and we can find $\theta \in \mbs^{n-1}$. We conclude that $\hat f(\zeta) = 0$ for $|\tau| \leq |\xi|$. By analyticity, $\hat f = 0$ and thus $f = 0. $
\epf

\begin{cor}\label{cor-uni}
Suppose $X_Mf = 0$ where $f$ satisfies the wave equation constraint \eqref{eq-cons} in which $f_1 \in H_{\comp}^{s+1}(\mbr^3), f_2\in H_{\comp}^{s}(\mbr^3), s\geq 0$ are compactly supported. Then $f = f_1 = f_2 = 0.$ 
\end{cor}
\bpf
Let $K = \supp f_1 \cup \supp f_2 \subset \mbr^3$. Let $I^+_c(K)$ be the chronological future of $K$ with respect to the Lorentzian metric induced by $c$. We know that there is a unique solution $f \in H^{s+1}(M)$ of \eqref{eq-cons}. By finite speed of propagation (or strong Huygens principle), the solution $f$ is supported in $I^+_c(K)\cap M$. Now we extend $f$ trivially to $\tilde f\in L^1_{\comp}(\mbr^4)$ and we regard $X_M$ as the light ray transform $X_{\mbr^4}$ on $\mbr^4$. We still have  $X_{\mbr^4} \tilde f = 0.$ By Theorem \ref{thm-inj}, we conclude that $f = 0$ on $\mbr^4$ so that $f = 0$ on $M$ and $f_1 = f_2 = 0$ on $\mcs_0.$
\epf

\bpf[Proof of Theorem \ref{thm-main1}]
The uniqueness part  is done in Corollary \ref{cor-uni}. So we prove the stability estimate below. We divide the proof into three steps. 

{\bf  Step 1: } Consider $c<1.$ From Proposition \ref{prop-para1}, we know that there are operators $A_1, A_2$ such that 
\beq
A_1 X_M f = f_1 + R_1 f_1 + R_1'f_2, \ \ A_2 X_M f = f_2 + R_2 f_1 + R_2'f_2
\eeq
and $R_i, R_i', i  =1, 2$ are all smoothing operators. We denote 
\beq
T
\begin{pmatrix}
f_1\\
f_2
\end{pmatrix}
 = \id 
 \begin{pmatrix}
f_1\\
f_2
\end{pmatrix}
+K
\begin{pmatrix}
f_1\\
f_2
\end{pmatrix}, \ \ K = \begin{pmatrix}
R_1 & R_1'\\
R_2 & R_2'
\end{pmatrix}.
\eeq
We consider  $T$ acting on $\mathcal{N}^s, s\geq 0.$ Then $K$ is compact from $\mathcal{N}^s$ to $\mathcal{N}^{s - \rho}, \rho \in \mbr$. So we have  the estimate 
\beq
 \|(f_1, f_2)\|_{\mathcal{N}^s} \leq \|A_1 X_M f\|_{H^{s+1}(\mbr^3)} + \|A_2 X_Mf\|_{H^{s}(\mbr^3)} +  C_\rho \|(f_1, f_2)\|_{\mathcal{N}^{s-\rho}} 
\eeq
for some constant $C_\rho$. 
Recall from Proposition \ref{prop-para1}  that
$A_1 = B_+ \circ (X_M\circ E_+)^* + B_- (X_M\circ E_-)^*$ and $A_2 = ic\lap^\ha (B_+\circ  (X_M\circ E_+)^* + B_- \circ  (X_M\circ E_-)^*)$.  Since the normal operator $(X_M E_\pm)^*X_M E_\pm$ are  pseudo-differential operators of order $-1$. By the $L^2$ estimate of pseudo-differential operators, we conclude that $X_M \circ E_\pm: H^s_{\comp}(\mbr^3)\rightarrow H^{s+\ha}_{\loc}(\mcc)$ is bounded.   Also, $(X_M \circ E_\pm)^*: H^{s}_{\comp}(\mcc) \rightarrow H^{s+\ha}_{\loc}(\mbr^3)$ is bounded. Therefore,   $A_1: H^{s+ \ha}_{\comp}(\mcc)\rightarrow H^{s}_{\loc}(\mbr^3)$ and $A_2: H^{s+ \ha}_{\comp}(\mcc)\rightarrow H^{s-\ha}_{\loc}(\mbr^3)$ are bounded.  
For $(f_1, f_2)\in \mathcal{N}^s$, we know from \eqref{eq-cauchysol} that $X_Mf = X_M E_+ h_1 + X_M E_- h_2$ and $h_1, h_2\in H^{s+1}(\mbr^3)$. Thus, $X_Mf \in H^{s+ 3/2}(\mcc)$ so we get  
\beqq\label{eq-estcomp}
 \|(f_1, f_2)\|_{\mathcal{N}^s} \leq C \|X_M f\|_{H^{s+3/2}(\mcc)} + C_\rho \|(f_1, f_2)\|_{\mathcal{N}^{s-\rho}} 
\eeqq
where $C_\rho > 0$ is a constant depending on $\rho.$  Note that the order is better than what claimed in the theorem for this case.

{\bf Step 2:} Consider $c = 1.$ It is convenient to work with $t_0>0$ which can be always arranged. For the Cauchy problem in Section \ref{sec-wave} with initial condition on $t = t_0$
\beqq\label{eq-scalar1}
\begin{gathered}
 \square   f   = 0, \quad \text{ on }  \mbr \times \mbr^3 \\
f = f_1, \quad \p_t f = f_2, \text{ on }   \{t_0\}\times \mbr^3,
\end{gathered}
\eeqq
it is known that  
\beq
U(t): (f_1, f_2) \rightarrow (f(t), \p_t f(t)), \quad t\in \mbr
\eeq
is bijective on $H^{s+1}(\mbr^3)\times H^s(\mbr^3)$. In fact, for \eqref{eq-scalar1}, $U(t)$ is a unitary operator with respect to the energy norm. 
We consider $(\tilde f_1, \tilde f_2) = U(-t_0)(f_1, f_2)$ which is the Cauchy data at $t = 0$ corresponding to $(f_1, f_2)$ at $t = t_0.$ 
Then we have 
\beqq\label{eq-energyest}
\begin{gathered}
C_1 (\|\tilde f_1\|_{H^{s+1}(\mbr^3)} + \|\tilde f_2\|_{H^{s}(\mbr^3)}) \leq \|f_1\|_{H^{s+1}(\mbr^3)} + \|f_2\|_{H^{s}(\mbr^3)} \\
 \leq C_2 (\|\tilde f_1\|_{H^{s+1}(\mbr^3)} + \|\tilde f_2\|_{H^{s}(\mbr^3)})
\end{gathered}
\eeqq
for some $C_1, C_2 > 0$, which follows from the energy estimate of the wave equation. We observe that the solution of \eqref{eq-scalar1} on $[t_0, t_1]\times \mbr^3$ can be expressed as 
\beq
f = \chi_{[t_0, t_1]} E(\tilde f_1, \tilde f_2) 
\eeq 
where $E(\tilde f_1, \tilde f_2) = E_+ \tilde h_1 + E_-\tilde h_2$ is the solution operator for the Cauchy problem from $t = 0$ in \eqref{eq-cauchysol} and $\tilde h_1, \tilde h_2$ correspond  to $\tilde f_1, \tilde f_2$, see Section \ref{sec-cauchy}. Therefore, we can apply Proposition \ref{prop-para2} to the operator $X_M\chi_{[t_0, t_1]}E_\pm$ with $t_0>0$. 

From Proposition \ref{prop-para2}, for any $\rho \in \mbn$,  there are operators $A_1, A_2$ such that 
\beq
A_1 X_M \chi_{[t_0, t_1]} f = \tilde h_1 + R_1 \tilde h_1 + R_1' \tilde h_2, \ \ A_2 X_M \chi_{[t_0, t_1]}  f = \tilde h_2 + R_2 \tilde h_1 + R_2' \tilde h_2
\eeq
and $R_i, R_i', i  =1, 2$ are FIOs of order $-\rho$. By the same argument in Step 1 and using Sobolev estimate of FIOs of canonical graph type, we have  
\beq
\begin{gathered}
 \|\tilde h_1\|_{H^s(\mbr^3)} +  \|\tilde h_2\|_{H^s(\mbr^3)} \leq \|A_1 X_M   f\|_{H^{s+1}(\mbr^3)} + \|A_2 X_M    f\|_{H^{s}(\mbr^3)} \\
 +  C_\rho (\|\tilde h_1\|_{H^{s-\rho}(\mbr^3)} +  \|\tilde h_2\|_{H^{s-\rho}(\mbr^3)})
 \end{gathered}
 \eeq
for some constant $C_\rho$.  Using \eqref{eq-fh}, we can change the estimate of $\tilde h_1, \tilde h_2$   to that of  $\tilde f_1, \tilde f_2$ and get 
\beq
\begin{gathered}
 \|(\tilde f_1, \tilde f_2)\|_{\mathcal{N}^s} \leq \|A_1 X_M \chi_{[t_0, t_1]}  f\|_{H^{s+1}(\mbr^3)} + \|A_2 X_M  \chi_{[t_0, t_1]}  f\|_{H^{s}(\mbr^3)} 
 +  C_\rho  \|(\tilde f_1, \tilde f_2)\|_{\mathcal{N}^{s-\rho}} 
 \end{gathered}
 \eeq
 Finally, using \eqref{eq-energyest}, we get 
 \beq
\begin{gathered}
 \|(f_1, f_2)\|_{\mathcal{N}^s} \leq \|A_1 X_M   f\|_{H^{s+1}(\mbr^3)} + \|A_2 X_M    f\|_{H^{s}(\mbr^3)} 
 +  C_\rho  \|(f_1, f_2)\|_{\mathcal{N}^{s-\rho}} 
 \end{gathered}
 \eeq
Now recall from the proof of Proposition \ref{prop-para2} that 
\beq\begin{gathered}
\begin{pmatrix}
A_1\\
A_2
\end{pmatrix}
 = W\begin{pmatrix}
\tilde A_1\\
\tilde A_2
\end{pmatrix}
\text{ where }
\tilde A_1 = Q^{1}_+ I^{\varphi_1} - Q^1_+ K^{\varphi_1, -} \tilde A_2, \quad 
\tilde A_2 = \sum_{k = 1}^3B_k(Q^{1}_+ I^{\varphi_1}     - Q^{2, k}_+ I^{\varphi_{2, k}} )
\end{gathered}
\eeq
in which $Q_+^1, Q_+^{2, k}\in \Psi^1(\mbr^3), k = 1, 2, 3$,  $B_k \in \Psi^0(\mbr^3), K^{\varphi_1, -}\in \Psi^{-1}(\mbr^3)$ and $W = \sum_{n = 0}^{\rho - 1}(-R)^n$ with elements of $R$ belonging to $I^{-1}(\mbr^3, \mbr^3; C_{t_0}) + I^{-1}(\mbr^3, \mbr^3; C_{t_1})+ I^{-1}(\mbr^3, \mbr^3; C_{-t_0})$ $+ I^{-1}(\mbr^3, \mbr^3; C_{-t_1})$. Using the estimate for pseudo-differential operators and FIOs of canonical graph type, we get 
\beq
\begin{gathered}
\|A_1 X_M    f\|_{H^{s+1}(\mbr^3)} + \|A_2 X_M    f\|_{H^{s}(\mbr^3)}\\
 \leq C \|I^{\varphi_1} X_M   f\|_{H^{s+2}(\mbr^3)} + C \sum_{k = 1}^3\|I^{\varphi_2, k} X_M    f\|_{H^{s+1}(\mbr^3)} \leq C\|X_M    f\|_{H^{s+2}}
\end{gathered}
\eeq
So in this case, we get 
\beqq\label{eq-estcomp1}
 \|(f_1, f_2)\|_{\mathcal{N}^s} \leq C \|X_M    f\|_{H^{s+2}(\mcc)} + C_\rho \|(f_1, f_2)\|_{\mathcal{N}^{s-\rho}}, \quad \rho \in \mbn.
\eeqq 

{\bf Step 3:} We get rid of the last term in \eqref{eq-estcomp} and \eqref{eq-estcomp1}. Let $\mck$ be a compact subset of $\mbr^3$ and denote by $\mathcal{N}^s(\mck)$ the function space consisting of $(f_1, f_2) \in \mathcal{N}^s$ supported in $\mck.$ Then the inclusion of $\mathcal{N}^s(\mck)$ into $\mathcal{N}^{s-\rho}(\mck), \rho>0$ is compact. We claim that 
\beq
 \|(f_1, f_2)\|_{\mathcal{N}^s(\mck)} \leq C \|X_M  f\|_{H^{s+2}(\mcc)}  
\eeq
for some $C >0$. We argue by contradiction. Assume the estimate without the error term is not true. We can get a sequence $(f^{(j)}_1, f^{(j)}_2), j = 1, 2, \cdots$ with unit norm in $\mathcal{N}^s(\mck)$ such that $X_M  f^{(j)}$ goes to 0 in $H^{s+2}(\mcc)$ as $j\rightarrow \infty$.  By \eqref{eq-estcomp} (for $(f_1, f_2)$ supported in $\mck$), we conclude that $1=\|(f_1^{(j)},f_2^{(j)})\|_{\mathcal{N}^s(\mck)}\leq C_\rho \|(f_1^{(j)},f_2^{(j)})\|_{\mathcal{N}^{s-\rho}(\mck)}$. This gives a weak limit $(f_1,f_2)$ in $\mathcal{N}^s(\mck)$ along a subsequence, which thus converges strongly in $\mathcal{N}^{s-\rho}(\mck)$. Therefore,  $\|(f_1, f_2)\|_{\mathcal{N}^{s-\rho}(\mck)}$ is bounded below by $1/C_\rho$, thus non-zero. However, $X_M  f = 0$ so $f = 0$ by the injectivity of $X_M$.  So $(f_1, f_2) = 0$ a contradiction. This finishes the proof. 
\epf

Finally, we prove  a stronger version of Theorem \ref{thm-main1} which allows lower order terms in the wave equation. We consider   differential operators of the form
\beq
P(x, t, D_x, \p_t) = \p_t^2 + c^2 \sum_{i = 1}^3D_{x_i}^2 + P_1(x, t, iD_x, \p_t) + P_0(x, t) 
\eeq
where $P_1$ is a first order differential operator with real valued smooth coefficients and $P_0$ is smooth.  Then we consider the Cauchy problem 
\beqq\label{eq-cons2}
\begin{gathered}
P(x, t, D_x, \p_t)f= 0 \ \ \text{ on } \intm\\
f = f_1, \ \ \p_t f = f_2, \ \ \text{ on } \mcs_0.
\end{gathered}
\eeqq
We remark that the equations for $\Phi$ in Section \ref{sec-wave} are  of this type. We   prove 
 \begin{theorem}\label{thm-main2}
Under the same assumptions as in Theorem \ref{thm-main1}, $X_M f$ uniquely determines $f$ and $f_1, f_2$ which satisfy \eqref{eq-cons2}. Moreover, there exists a $C> 0$ such that 
\beq
\|(f_1, f_2)\|_{\mathcal{N}^s} \leq C  \|X_M   f\|_{H^{s+ 2}(\mcc)} \text{ and } \|f\|_{H^{s+1}(M)} \leq C  \|X_M    f\|_{H^{s+ 2}(\mcc)} 
\eeq
where  $\mcc$ is the set of light rays on $M$.
 \end{theorem}
 \bpf
 The proof follows the same arguments as for Theorem
 \ref{thm-main1}. So  we just point out what needs to be
 modified. When the wave equation contains lower order terms, one can
 construct parametrices $E_\pm$ for the Cauchy problem, see
 \cite[Chapter 5]{Du}. These are Fourier integral operators and can be
 represented by oscillatory integrals. So the construction in Section
 \ref{sec-cauchy} works through, and the analysis for 
 $X_ME_\pm$ is the same as the standard wave equation case. However,
 we do need to justify the ellipticity of the involved operators in
 Lemma \ref{lm-elliptic1} and Proposition  \ref{prop-compose}. 
 We remark that ellipticity of the solution itself
 is standard, and follows simply from the principal symbol satisfying
 a transport equation. 
We follow the parametrix construction in Tr\`eves \cite[Section 1, Chapter VI]{Tre} to check this in a transparent manner. 

We look for operators $E_j, j = 0, 1$ such that 
\beq 
\begin{gathered}
P(x, t, D_x, \p_t)E_j  = 0 \ \ \text{ on } \intm\\
\p_t^k E_j  = \delta_{kj},  k = 0, 1,  \ \ \text{ on } \mcs_0 .
\end{gathered}
\eeq 
Here, for $j = 0, 1$ we have 
\beq
E_j  f(x)  = (2\pi)^{-3} \int_{\mbr^3} e^{i\phi_0(x, t, \xi)} a_{j0}(x, t, \xi) \hat f(\xi) d\xi  + (2\pi)^{-3} \int_{\mbr^3} e^{i\phi_1(x, t, \xi)} a_{j1}(x, t, \xi) \hat f(\xi) d\xi+ R_j(t)f(x)
\eeq
where $R_j$ are smoothing operators, see \cite[(1.37)]{Tre}. The phase functions are 
\beq
\phi_0(x, t, \xi) = x\cdot \xi + ct |\xi|, \ \  \phi_1(x, t, \xi) = x\cdot \xi - ct|\xi|.
\eeq
The amplitude can be written as $a_{jk}(x, t, \xi) = \sum_{l = 0}^\infty a_{jkl}(x, t, \xi)$ and each $a_{jkl}$ is homogeneous of degree $-j-l$ for $|\xi|$ large. 
Before we look into the structures that we need of the amplitude, we  find the initial values of the leading order term $a_{jk0}$ at $t = t_0$. They satisfy (see \cite[(1.53)]{Tre})
\beq
a_{000}(x, t, \xi) = \ha, \ \ a_{010}(x, t, \xi) =  \ha, \ \ a_{100}(x, t, \xi) = \frac{1}{2ic|\xi|}, \ \ a_{110}(x, t, \xi) = -\frac{1}{2ic|\xi|}.
\eeq 
The amplitudes satisfy first order equations which are deduced from (see \cite[(1.39)]{Tre})
\beq
P(x, t, D_x + \p_x \phi_k, \p_t + i\p_t \phi_k)a_{jk}(x, t, \xi) = 0.
\eeq
For the leading order term, we get 
\beqq\label{eq-transport}
\p_\tau P_2(x, t, \p_x \phi_k, i\p_t \phi_k) \p_t a_{jk0} + \sum_{\nu = 1}^3 \p_{\xi_\nu}P_2(x, t, \p_x\phi_k, i\p_t\phi_k) D_{x^\nu} a_{jk0} + C(\phi_k; x, t, \xi)a_{jk0} = 0
\eeqq
and the $C$ term in this case is (the sub-principal symbol of $P$)
\beq
C(\phi_k; x, t, \xi) = P_1(x, t, i\p_x \phi_k, i\p_t \phi_k). 
\eeq
Note that $P_1$ has real valued coefficients and is homogeneous of degree one in $i\p_x\phi_k, i\p_t\phi_k$. Dividing by $i = \sqrt{-1}$, we see that equation \eqref{eq-transport} is a first order linear equation with real valued coefficients. Solving the equation amounts to solving a ODE along the integral curve and the solution $a_{jk0}$ will be  positive scalar multiples of the initial conditions hence not only non-vanishing, but is real or purely imaginary depending on its initial value.

Finally, we can represent the solution to \eqref{eq-cons2} as
\beq
f(x, t) = E_0 f_1 + E_1 f_2 = E_{+} h_1 + E_-{h_2} 
\eeq
where 
\beqq\label{eq-epm}
\begin{gathered}
E_{+}h = (2\pi)^{-3} \int_{\mbr^3} e^{i(x\cdot \xi + ct|\xi|)} (a_{00}(x, t, \xi) + 2ic|\xi| a_{10}(x, t, \xi)) \hat h(\xi) d\xi \\
  = (2\pi)^{-3} \int_{\mbr^3} e^{i(x\cdot \xi + ct|\xi|)}  a_+(x, t, \xi)  \hat h(\xi) d\xi  \\
E_{-}h = (2\pi)^{-3} \int_{\mbr^3} e^{i(x\cdot \xi - ct|\xi|)} (a_{01}(x, t, \xi) - 2ic|\xi| a_{11}(x, t, \xi)) \hat h(\xi) d\xi \\
= (2\pi)^{-3} \int_{\mbr^3} e^{i(x\cdot \xi - ct|\xi|)}  a_-(x, t, \xi)  \hat h(\xi) d\xi 
\end{gathered}
\eeqq
and 
\beq
h_1 = f_1 + \frac{1}{2ic}\lap^{-\ha}   f_2, \ \ h_2 = f_1 -   \frac{1}{2ic}\lap^{-\ha} f_2.
\eeq
We see that the leading order terms of $a_+, a_-$ are all positive. From these oscillatory integral representations, it is easy to see that Lemma \ref{lm-elliptic1} holds for $c<1$. For Proposition \ref{prop-compose}, we see that the principal symbol of $k_+$ is given by
\beq
k_{+, -1}(x, \xi) = 2\pi i |\xi|^{-1}  \varphi(-\xi/|\xi|) \int_0^{t_1}t^{-1} \chi_\eps(t) a_{+, 0}(t, x)dt  
\eeq 
 where $a_{+, 0}$ is the in the expansion $a_+\sim \sum_{k = 0}^\infty a_{+, k}(t, \xi) |\xi|^{-1-k}$. So $k_{+, -1}(x, \xi)$ is non-vanishing. Thus the operator $I^\varphi X_M\chi_\eps E_+$ is elliptic. 
The rest of the proof is the same as for Theorem \ref{thm-main1}.   
 \epf

 \section{Small perturbations of the Minkowski spacetime}\label{sec-pert}
We consider metric perturbations $g_\delta = g_M + h$ with $h = \sum_{i, j = 0}^3 h_{ij} dx^idx^j$. We assume that 
\begin{enumerate}
\item[(A1)] $h$ is a symmetric two tensor smooth on $M$; 
\item[(A2)]  for $\delta > 0$ small,  the seminorm  $\|h_{ij}\|_{C^{3}} $ $= \sup_{(t, x)\in M}\sum_{|\alpha|\leq 3}|\p^\alpha h_{ij}(t, x)|  <  \delta, i, j = 0, 1, 2, 3.$
\end{enumerate}
 Without loss of generality, we can assume that $h$ is extended to some larger manifold $\tilde M = (\tilde t_0, \tilde t_1)\times \mbr^3$ such that $M\subset \tilde M$ and (A2) holds on $\tilde M$. 
In this section, we study the inverse problem on $(M, g_\delta)$ for $\delta$ sufficiently small. Note that in this case, light rays may not follow straight lines and the injectivity of the light ray transform on scalar functions is not known. We will show that by using a perturbation argument on the Fourier integral operator level, one can obtain the same determination result as for the Minkowski case.  

We start with the light-like geodesics on $(M, g_\delta)$ and their parametrizations. Let $\gamma(s)$ denote a light-like geodesic from $\mcs_0$. It satisfies  
\beqq\label{eq-geo}
\begin{gathered}
\p_s^2 \gamma^k(s) + \Gamma^{k}_{ij} \p_s\gamma^i(s) \p_s\gamma^j(s) = 0\\
\gamma(0) = (0, y), \p_s \gamma(0) = (\beta, v)
\end{gathered}
\eeqq
where $\Gamma^{k}_{ij}$ is the Christoffel symbol for $g_\delta$, $v\in \mbs^2$ and $\beta$ is such that  $g_\delta(\beta, v) = 0$ and $(\beta, v)$ future pointing. It is known, see for example \cite{AM},  that \eqref{eq-geo} is equivalent to a first order system on $T^*M.$ Here, $M$ is regarded as a submanifold of $\tilde M. $
We use $(t, x)$ and $(\tau, \xi)$ for the local coordinates on $T^*M.$ Consider the Hamiltonian 
\beq
p(t, x, \tau, \xi) = \ha g^*_\delta(\tau, \xi) = \ha g_M^*(\tau, \xi) + H(t, x, \tau, \xi) = \ha(-|\tau|^2 + \sum_{i = 1}^3|\xi_i|^2) + H(t, x, \tau, \xi). 
\eeq 
Here, $2H$ is the perturbation of the dual metric  corresponding to the perturbation $h$. 
Let $\Xi = (\tau, \xi)$, then $H(t, x, \Xi) = \sum_{i, j = 0, 1, 2, 3} H_{ij}(t, x) \Xi_i\Xi_j$ is homogeneous of degree two in $\Xi$ and the seminorm $\|H_{ij}\|_{C^{3}} < C\delta$ for some constants $C$. We denote the Hamilton vector field by $H_p$. Let $(t(s), x(s), \tau(s), \xi(s))$ be an integral curve of $H_p$ in the characteristic set $\Sigma_p = \{(t, x, \tau, \xi) \in T^*M:  p(t, x, \tau, \xi) = 0\}$, called  null-bicharacteristics.   With $\gamma(s) = (t(s), x(s))$,  \eqref{eq-geo} can be converted to 
\beqq\label{eq-bichar}
\begin{gathered}
\frac{d t}{ds}  = \frac{\p p}{\p \tau} = - \tau + \p_\tau H(t, x, \tau, \xi); \quad \frac{d x_i}{ds} = \frac{\p p}{\p \xi_i} =  \xi_i + \p_{\xi_i} H(t, x, \tau, \xi)\\
\frac{d \tau}{ds}  = -\p_{t} H(t, x, \tau, \xi); \quad \frac{d \xi_i}{ds} = -\p_{x_i}H(t, x, \tau, \xi), \quad i = 1, 2, 3\\
t(0) = t_0= 0, \quad x_i(0) = y_i, \quad \tau(0) = \tau_0, \quad \xi_i(0) = \xi_{0, i}.
\end{gathered}
\eeqq
Here, $(\tau_0, \xi_0)$ is the cotangent vector obtained from $(\beta, v)$ using $g_\delta$ and we also denote it by $(\tau_0, \xi_0) = (\beta, v)^\flat$.  If we consider the system for the Minkowski metric namely $H = 0$, then  $\beta = 1$ and the covector $(\tau_0, \xi_0) = (-1, v)$. \eqref{eq-bichar} becomes 
\beqq\label{eq-bicharm}
\begin{gathered}
\frac{dt}{ds}  =  - \tau,  \quad \frac{d x_i}{ds} = \xi_i,  \quad  \frac{d \tau }{ds} = 0, \quad \frac{d \xi_i}{ds} = 0, \quad i = 1, 2, 3\\
t(0) = 0, \quad x_i(0) = y_i, \quad \tau(0) =  -1, \quad \xi_i(0) = v_i. 
\end{gathered}
\eeqq
We see that $x(s) = (s, y + s v), t(s) = s$, which agrees with our parametrization used previously. Now we have the following result. 
\begin{lemma}\label{lm-bichar}
For $\delta > 0$ sufficiently small, the set of light rays on $(M, g_\delta)$ is given by $\mcc_\delta = \{\gamma = (t, x(t, y, v)):  (y, v)\in \mcs_0 \times \mbs^2, t\in [t_0, t_1]\}$, where $x$ is a smooth function of $t, y, v$.  Moreover,  we have 
\beq
\|x(t, y, v) - (y + tv) \|_{C^2} < C \delta 
\eeq
for some constant $C$. 
\end{lemma}
\bpf
For $v\in \mbs^2$, the co-vectors $(\tau_0, \xi_0) = (\beta, v)^\flat$ are in a bounded set of $\mbr^4$. We assume that $|(\tau_0,  \xi_0)| < M_1$. We also notice that $\tau_0$ is away from zero, say $|\tau_0| > M_0 > 0$. Then we consider $(\tau, \xi)$ such that $|(\tau, \xi) - (\tau_0, \xi_0)| < M_0/2$ so that $|(\tau, \xi)|< M\doteq M_1 + M_0/2$ and $|\tau| > M_0/2$. Consider the system \eqref{eq-bichar}. Because $H$ is homogeneous of degree two in $(\tau, \xi)$, for $|(\tau, \xi)|< M$ and for $\delta > 0$ sufficiently small, we see that $\frac{dt}{ds} \neq 0$. Therefore, we can take $t$ as the parameter and convert \eqref{eq-bichar} to 
\beqq\label{eq-bichar2}
\begin{gathered}
\frac{d s}{dt}  =\frac{1}{- \tau + \p_\tau H(t, x, \tau, \xi)}; \quad \frac{d x_i}{dt} = \frac{\xi_i + \p_{\xi_i} H(t, x, \tau, \xi)}{- \tau + \p_\tau H(t, x, \tau, \xi)}\\
\frac{d \tau}{dt}  = \frac{-\p_{t} H(t, x, \tau, \xi)}{- \tau + \p_\tau H(t, x, \tau, \xi)}; \quad \frac{d \xi_i}{dt} = \frac{-\p_{x_i}H(t, x, \tau, \xi)}{- \tau + \p_\tau H(t, x, \tau, \xi)}, \quad i = 1, 2, 3\\
s(0) = 0, \quad x_i(0) = y_i, \quad \tau(0) = \tau_0, \quad \xi_i(0) = \xi_{0, i}.
\end{gathered}
\eeqq
The system corresponding to \eqref{eq-bicharm} is 
\beqq\label{eq-bichar3} 
\begin{gathered}
\frac{d s}{dt}  =\frac{1}{- \tau}; \quad \frac{d x_i}{dt} = \frac{\xi_i }{- \tau }, \quad 
\frac{d \tau}{dt}  = 0; \quad \frac{d \xi_i}{dt} = 0 , \quad i = 1, 2, 3\\
s(0) = 0, \quad x_i(0) = y_i, \quad \tau(0) = -1, \quad \xi_i(0) = v.
\end{gathered}
\eeqq 

Let $(\tilde t, \tilde x, \tilde \tau, \tilde \xi)$ be the solution of  \eqref{eq-bichar3} and $(t, x, \tau, \xi)$ satisfy \eqref{eq-bichar2}. Then let $u = (t - \tilde t, x - \tilde x, \tau - \tilde \tau, \xi - \tilde \xi)$. We see that $u$ satisfies the system 
\beq
\begin{gathered}
\frac{d u}{ds} = F(u)\\
u(0) =   u_0,
\end{gathered}
\eeq
where $F$ is smooth and $|F(u)| < C\delta$, $|u_0|< C\delta$ for generic constant $C$. Now it follows from standard ODE theorems, see for instance \cite[Theorem 1.2.3]{Ho5} that for $\delta$ sufficiently small, there is a unique $C^\infty$ solution $u$ on $[t_0, t_1]$ and $|u| \leq C\delta.$ Higher order estimates can be obtained similarly. This finishes the proof.  
\epf

Now we consider the light ray transform $X_\delta$ on $(M, g_\delta)$. The parametrization of the light rays is not unique, although all choices give rise to equivalent analysis for our purpose. Perhaps the most natural parameterization is to use the cosphere bundle on $\mcs_0$ of the induced metric. Let $\bar g_\delta$ be the induced Riemannian metric of $g_\delta$ on $\mcs_0$. For $y\in \mcs_0,$ let $\mbs_{\delta, y}^2 = \{v \in T\mcs_0: \bar g_\delta(v, v) = 1\}$. For $v\in \mbs_{\delta, y}^2$, there is a   unique future pointing light-like vector $(v_0, v)$ at $y$. In particular, $v_0$ is close to $1$ for $\delta$ small.  Then the light ray from $(0, y)$ in direction $(v_0, v)$ is parametrized by $\gamma_{y, v}(s) = \exp_{(0, y)} s(v_0, v), s\in [0, s_1]$
 where $s$ is the affine parameter such that $\gamma_{y, v}(0)  = (0, y)\in \mcs_0$ and $\gamma_{y, v}(s_1)\in \mcs_1$. In this parametrization, we can write 
\beqq\label{eq-xraypert}
\begin{gathered}
X_{\delta}f(y, v) = \int_0^{s_1} f(\gamma_{y, v}(s))d s. 
\end{gathered}
\eeqq 

Now we can identify $\mbs_{\delta, y}^2$ with $\mbs_{y}^2$ via a diffeomorphism. By the above Lemma \ref{lm-bichar}, $s$ is a smooth function of $y, t$ and $v\in \mbs^2$ so  we can use $t$ variable to parametrize the light rays. We have 
\beq 
\begin{gathered}
X_{\delta}f(y, v) = \int_0^{t_1} w(y, v, t) f(t, x(t, y, v))dt, \quad y\in \mcs_0, v\in \mbs^2, 
\end{gathered}
\eeq 
where $w$ is a weight coming from the change of variables. In fact, $w$ is smooth and close to $1$ for $\delta$ sufficiently small. $w$ only mildly affects the argument, changing the elliptic principal symbol of the final operator $X_\delta \circ E_+$ in \eqref{eq-XME1}, thus maintaining ellipticity. For simplicity, we will ignore it in the follows and take 
\beqq 
\begin{gathered}
X_{\delta}f(y, v) = \int_0^{t_1}  f(t, x(t, y, v))d t
 = (2\pi)^{-3}\int_{\mbr^3} \int_{\mbr^3} \int_{0}^{t_1} e^{i ( (x(t, y, v) - z)\cdot \eta)} f(t, z)  dt dz  d\eta.
\end{gathered}
\eeqq 
This is the parametrization of $X_\delta$ we work with in the rest of this section. 
The Schwartz kernel of $X_{\delta}$ is  the delta distribution on $\mcc \times \intm$ supported on the point-line relation $Z_\delta$ defined by  
\beq
\begin{gathered}
Z_\delta = \{(\gamma, q)\in \mcc\times \intm: q\in \gamma\} 
= \{(y, v, (t, x)) \in \mbr^3 \times \mbs^2 \times \intm:  x = x(t, y, v)\}. 
\end{gathered}
\eeq

Next, let $\square_{g_\delta}$ be the d'Alembert operator  on $(M, g_\delta)$ and  we consider the second order operator
\beqq
P_\delta(x, t, D_x, \p_t) = \square_{g_\delta} + P_1(x, t, iD_x, \p_t) + P_0(x, t) 
\eeqq
where $P_1$ is a first order differential operator with real valued smooth coefficients and $P_0$ is smooth.  Then we consider the Cauchy problem 
\beqq\label{eq-cons3}
\begin{gathered}
P_\delta(x, t, D_x, \p_t)f= 0 \ \ \text{ on } \intm \\
f = f_1, \ \ \p_t f = f_2, \ \ \text{ on } \mcs_0.
\end{gathered}
\eeqq
We remark that for sufficiently small metric perturbations, the operators $\square_{g_\delta}$ and $P_\delta$ are both strictly hyperbolic with respect to $\mcs_0.$ Therefore, as in previous sections, the parametrix construction of Duistermaat-H\"ormander can be applied. In general, the parametrix does not have a global oscillatory integral representation on $M$. However, we show below that  for sufficiently small perturbations of the Minkowski spacetime, this is possible.  

The parametrix construction is the same as in the previous section.  We look for operators $E_j, j = 0, 1$ such that 
\beq 
\begin{gathered}
P_\delta (x, t, D_x, \p_t)E_j  = 0 \ \ \text{ on } \intm \\
\p_t^k E_j  = \delta_{kj},  k = 0, 1,  \ \ \text{ on } \mcs_0 .
\end{gathered}
\eeq 
For $j = 0, 1$ we have 
\beq
E_j  f(x)  = (2\pi)^{-3} \int_{\mbr^3} e^{i\phi_+(x, t, \xi)} a_{j,+}(x, t, \xi) \hat f(\xi) d\xi  + (2\pi)^{-3} \int_{\mbr^3} e^{i\phi_-(x, t, \xi)} a_{j, -}(x, t, \xi) \hat f(\xi) d\xi+ R_j(t)f(x)
\eeq
where $R_j$ are smoothing operators, see \cite[(1.37)]{Tre}.  We follow  Tr\`eves \cite{Tre} to find the phase functions $\phi(t, x, \xi)$ for $(t, x)\in (t_0, t_1) \times \mbr^3$, $\xi \in \mbr^3$.  The phase function should satisfy the eikonal equation 
\beq
\begin{gathered}
p(\nabla \phi) =  -|\p_t\phi|^2 + |\p_x \phi|^2 + H(\p_t \phi, \p_x \phi) = 0.
 \end{gathered}
 \eeq
By the strict hyperbolicity, there are two solutions for $\p_t\phi$ denoted by $\p_t \phi = \la_\pm(t, x, \p_x \phi)$ and $\la_\pm$ are smooth functions and homogeneous of degree one in $\p_x \phi$. We take initial conditions $\p_t \phi = x\cdot \xi, \xi\in \mbr^3$ at $t = 0$. Below, we consider $\la_+$. The treatment for $\la_-$ is identical. We consider the Hamilton-Jacobi equation
 \beqq\label{eq-HJ}
 \begin{gathered}
 \frac{dx}{dt}  = - \p_\eta \la_{+}(t, x, \eta), \quad \frac{d \eta}{dt}  =  \p_x \la_+(t, x, \eta)\\
 x(0) = y, \quad \eta(0) = \xi, \quad y \in \mbr^3, \xi\in \mbr^3\backslash 0. 
 \end{gathered}
 \eeqq
We denote the solution by $x(t, y, \xi), \xi(t, y, \xi)$.  Then  the phase function is 
\beqq\label{eq-phase}
\phi_+(t, x, \xi) = x\cdot \xi + \int_0^t \la_+(s, x, \eta(s, y, \xi))ds. 
\eeqq
 Here, one can express $y$ in terms of $x$, see \cite[Section 2, Chapter VI]{Tre} for more details.  
 For the Minkowski spacetime, we know $\la_+ =  -|\xi|$ so that \eqref{eq-HJ} becomes
 \beqq\label{eq-HJm}
 \begin{gathered}
 \frac{dx}{dt}  =   \xi/|\xi|, \quad \frac{d \eta}{dt}  =  0\\
 x(0) = y, \quad \eta(0) = \xi. 
 \end{gathered}
 \eeqq
The solution is simply $x(t) = y  +  t\xi/|\xi|, \eta(t) = \xi$ and the phase function is $\phi_0(t, x, \xi) = x\cdot \xi + t|\xi|$. 
Using the same argument as for Lemma \ref{lm-bichar}, we get
\begin{lemma}\label{lm-phaseest}
For $\delta > 0$ sufficiently small, there is a unique smooth solution $(x(t, y, \xi), \eta(t, y, \xi))$ to \eqref{eq-HJ} for $t\in[t_0, t_1], y\in \mbr^3, \xi \in \mbr^3\backslash 0$, and they satisfy 
\beq
\| x(t, y, \xi) - (y - t \xi/|\xi|)  \|_{C^2} < C\delta,  \quad \|\eta(t, y, \xi)/|\xi| - \xi/|\xi|\|_{C^2} < C \delta 
\eeq
for some constant $C > 0.$ It follows that the phase function $\phi_+$ in \eqref{eq-phase} is also smooth and satisfies 
\beq
\|\phi_+(t, x, \xi) - (x\cdot \xi + t|\xi|)\|_{C^2} < C\delta|\xi|.
\eeq
\end{lemma}

We remark that similar argument was used in \cite{SU1} for a backscattering problem. Using this lemma, we can represent the solution to \eqref{eq-cons3} as
\beq
f(x, t) = E_0 f_1 + E_1 f_2 = E_{+} h_1 + E_-{h_2} 
\eeq
where 
\beqq\label{eq-epm2}
\begin{gathered}
E_{+}h  
  = (2\pi)^{-3} \int_{\mbr^3} e^{i  \phi_+(t, x, \xi)  }  a_+(x, t, \xi)  \hat h(\xi) d\xi   \\
E_{-}h  
= (2\pi)^{-3} \int_{\mbr^3} e^{i \phi_-(t, x, \xi)}  a_-(x, t, \xi)  \hat h(\xi) d\xi.  
\end{gathered}
\eeqq
The $a_\pm$ and $h_1, h_2$ are the same as in \eqref{eq-epm}. 
\newline

With these preparations, we now  state and prove our  main result in this section. 
 \begin{theorem}\label{thm-main3}
Consider $(M, g_\delta)$ which satisfy the assumptions (A1), (A2) in the beginning of this section. Assume that $(f_1, f_2)\in \mathcal{N}^s, s \geq 0 $, and $f_1, f_2$ are supported in a compact set $\mck$ of $\mcs_0$. For $\delta \geq 0$ sufficiently small, $X_\delta f$ uniquely determines $f$ and $f_1, f_2$ which satisfy \eqref{eq-cons3}. Moreover, there exists   $C> 0$ such that 
\beq
\|(f_1, f_2)\|_{\mathcal{N}^s} \leq C  \|X_\delta    f\|_{H^{s+ 2}(\mcc_\delta)} \text{ and } \|f\|_{H^{s+1}(M)} \leq C  \|X_\delta     f\|_{H^{s+ 2}(\mcc_\delta)} 
\eeq
where  $\mcc_\delta$ is the set of light rays on $(M, g_\delta)$.
 \end{theorem}
 
\bpf
We examine the arguments in Section \ref{sec-local2} and Section \ref{sec-full} and point out what needs to be modified. We consider the composition of $X_\delta$ and $E_+$ defined in \eqref{eq-epm2}. We have 
\beq
\begin{gathered}
X_{\delta}f(y, v)  
 = (2\pi)^{-3}\int_{\mbr^3} \int_{\mbr^3} \int_{0}^{t_1} e^{i ( (x(t, y, v) - x')\cdot \eta)} f(t, x')  dt dx'  d\eta
 \end{gathered}
\eeq 
and  
\beq
E_+(f)(t, x') =(2\pi)^{-3} \int_{\mbr^3}\int_{\mbr^3} e^{i(\phi_+(t, x', \xi) - z \cdot \xi)} a_+(t, x', \xi) f(z) dz d\xi. 
\eeq
Consider the integral operator $I^\varphi$ defined in Section \ref{sec-local2}. Using the oscillatory integral representations, we have 
\beqq\label{eq-XME1}
 \begin{gathered}
I^\varphi X_{\delta}  \chi_\eps E_+ f(y, v)   \\
= (2\pi)^{-6} \int_{\mbs^2}\int_{\mbr^3}\int_{\mbr^3}\int_0^{t_1}\int_{\mbr^3}\int_{\mbr^3} e^{i (x(t, y, v)\cdot \eta - x'\cdot \eta+ \phi_+(t, x', \xi) - z\cdot \xi)}  \varphi(v) a(t, x', \xi) \chi_\eps(t) f(z) dz d\xi dt dx' d\eta dv \\
\end{gathered}
\eeqq
We  write the phase function as $\Phi = \phi + \psi$ in which  
\beq
\phi(y, z; \xi, \eta, x', t, v) = (y - x')\cdot \eta + t v\cdot \eta+ (x' -z)\cdot \xi +   t |\xi| 
\eeq 
and $\psi$ is a smooth function and homogeneous of degree one in $\xi, \eta$. In particular, $\Phi$ is a small perturbation of $\phi$. As in Proposition \ref{prop-compose}, we first consider the integration in $x', \eta, v$ in \eqref{eq-XME1}. As shown in Proposition \ref{prop-compose}, the phase function $\phi$ in these variables is non-degenerate. Since $\psi$ is a small perturbation of $\phi$, for $\delta$ sufficiently small, we see that $\Phi$ in $x', \eta, v$ variables is also non-degenerate. Note that 
\beqq\label{eq-station}
\p_{x'}\Phi = -\eta + \p_{x'}\phi_+(t, x', \xi), \quad \p_{\eta} \Phi = x(t, y, v) - x', \quad \p_{v}\Phi = \p_vx(t, y, v)\cdot \eta 
\eeqq
For the stationary points, we see that $x' = x(t, y, v)$ so $(t, x)$ is on the light ray from $(0, y)$ in direction $(1, v)$. 
Let $\tau$ satisfy $p_\delta(t, x, \tau, \eta) = 0$. From $\eta = \p_{x'}\phi_{+}(t, x', \xi)$   we see that $(t, x', \tau, \eta)$ is on the bicharactersitics from $(y, \xi)$. Since there is no conjugate points, we get $v = \pm\xi/|\xi|$.  
Thus at the stationary points, the phase function becomes
\beq
\Phi(y, z, t, \xi) = \phi_+(t, x(t, y, \pm \xi/|\xi|), \xi) - z\cdot \xi
\eeq
After integrating in $x', \eta, v$,  the Schwartz kernel becomes
\beqq\label{eq-ivint} 
 \begin{gathered}
I^\varphi X_{\delta}  \chi_\eps E_+ (y, z)    
= (2\pi)^{-3}  \int_0^{t_1}\int_{\mbr^3}  e^{i (\phi_+(t, x(t, y,  \xi/|\xi|), \xi) - z\cdot \xi)} k^\delta_+(t, \xi) d\xi dt  \\ 
+ (2\pi)^{-3}  \int_0^{t_1}\int_{\mbr^3} e^{i (\phi_+(t, x(t, y, - \xi/|\xi|), \xi) - z\cdot \xi)} k^\delta_-(t, \xi) d\xi dt  
\end{gathered}
\eeqq
where $k^\delta_\pm$ are small perturbations of $k_\pm$ in \eqref{eq-kplus} and \eqref{eq-kminus} of Proposition \ref{prop-compose}. Finally, we integrate in $t$. For the second integral in \eqref{eq-ivint}, the phase function is a small perturbation of 
\beq
(y - z) \cdot \xi + 2t|\xi|
\eeq
thus as in Proposition \ref{prop-compose}, the integral is $O(|\xi|^{-\infty})$. For the first integral of \eqref{eq-ivint}, we need to examine the phase function at the stationary points. Using \eqref{eq-phase}, we get 
\beq  
\Phi(y, z, t, \xi) = x\cdot \xi + \int_0^t \la_+(s, x, \xi)ds - z\cdot \xi
\eeq 
where $x = x(t, y, \xi/|\xi|)$. Taking $\xi$ derivative, we get 
\beq
\begin{gathered}
\p_\xi\Phi(y, z, t, \xi) = (x - z)  + \p_\xi x \cdot \xi  + \int_0^t (\p_\eta \la_+(s, x(s, y, \xi/|\xi|), \xi)\p_\xi x  + \p_\eta \la_+) ds \\
 = (x - z) + \p_\xi x \cdot \xi  + \int_0^t  -\frac{dx}{ds}(s, x, \xi) ds
\end{gathered}
\eeq
where we used the stationary point condition \eqref{eq-station} and $\p_\eta \la_+ = -dx/dt$. Note that $\p_\xi x$ is the Jacobi field, and because $x(t, y, \xi/|\xi|)$ is a light-like geodesic, $\p_\xi x\cdot \xi = 0$, see Lemma 3.1 and Lemma 3.4 of \cite{LOSU1}. Therefore,  $
\p_\xi\Phi(y, z, t, \xi)  
 = (x - z)$ and $\Phi(y, z, t, \xi) = (y - z)\cdot \xi + \tilde \Phi(y, z, t)$ where $\tilde \Phi$ is small. Finally, integrating in $t$ of the first integral of \eqref{eq-ivint} gives a pseudo-differential operator of order $-1$ and the principal symbol $k^\delta_{+, -1}$ is a small perturbation of $k_{+, -1}(\xi)$ in Proposition \ref{prop-compose}.  This implies that Proposition \ref{prop-compose} hold for the small perturbations. 
 
To see that the analogous result of Proposition \ref{prop-para2} holds for small perturbations, it suffices to examine the kernel \eqref{eq-ivint} in which $\chi_\eps$ is replaced by $\chi_{[\eps, t_1]}$ 
 \beqq\label{eq-ivint1} 
 \begin{gathered}
I^\varphi X_{\delta}  \chi_{[\eps, t_1]} E_+ (y, z)    
= (2\pi)^{-3}  \int_\eps^{t_1}\int_{\mbr^3} e^{i (\phi_+(t, x(t, y,  \xi/|\xi|), \xi) - z\cdot \xi)} k^\delta_+(t, \xi) d\xi dt  \\ 
+ (2\pi)^{-3}  \int_\eps^{t_1}\int_{\mbr^3} e^{i (\phi_+(t, x(t, y, - \xi/|\xi|), \xi) - z\cdot \xi)} k^\delta_-(t, \xi) d\xi dt  
\end{gathered}
\eeqq
 The first integral still gives a pseudo-differential operator as shown above. For the second integral, integration by parts in $t$ gives an oscillatory integral of the form 
 \beqq\label{eq-ivint2}
 \int_{\mbr^3}  e^{i (\phi_+(\eps, x(\eps, y, - \xi/|\xi|), \xi) - z\cdot \xi)} a(\xi) d\xi  +  \int_{\mbr^3}  e^{i (\phi_+(t_1, x(t_1, y, - \xi/|\xi|), \xi) - z\cdot \xi)} b(\xi) d\xi
 \eeqq
 where $a, b$ are symbols of order $-2$. Here, we used that $\phi_+$  is homogeneous of degree one in $\xi.$ To see that these are FIOs of canonical graph type, we use the characterization in \cite[page 26]{Ho4} which says that an oscillatory integral with phase $\phi(x, \eta) - x \cdot \eta$  is an FIO whose canonical relation is a canonical graph if and only if $\det \frac{\p^2\phi}{\p x\p \eta} \neq 0.$ Since $\phi_+(\eps, x(\eps, y, - \xi/|\xi|), \xi)$ is a small perturbation of $y\cdot \xi  + 2 \eps |\xi|$ and $\det \frac{\p^2 }{\p y\p \xi} (y\cdot \xi  + 2 \eps |\xi|) = -1 \neq 0$, we conclude that for $\delta$ sufficiently small, the first integral in \eqref{eq-ivint2} gives an FIO of canonical graph type. The same is true for the second integral. Thus Proposition \ref{prop-para2} holds for small perturbations.

Now, the proof of Theorem \ref{thm-main1} in Section \ref{sec-full} go through line by line, except the injectivity  of $X_\delta$. In particular, we have the  estimate as  \eqref{eq-estcomp}
\beq 
 \|(f_1, f_2)\|_{\mathcal{N}^s} \leq C \|X_\delta   f\|_{H^{s+2}(\mcc_\delta)} + C_\rho \|(f_1, f_2)\|_{\mathcal{N}^{s-\rho}} 
\eeq
where $C_\rho$ is a constant depending on $\rho.$ To get rid of the last term, we use the following argument, see \cite[Section 2.7]{Va}. 
Notice that given $s,\rho$  and for some fixed small $\delta_0$, if we consider all metric $g$ such that  $\|g - g_M\|_{C^{3}} \leq \delta_0$ , then the above estimate is uniform (a fixed constant $C_\rho$ works for all such metrics) by the uniformity of the construction.  
Now suppose there is no $\delta$ such that for all metrics within $\delta$ of the Minkowski metric $g_M$ (in the Fr\'echet space sense) the transform is injective. Let $F^j=(f^j_1,f^j_2), j = 1, 2, \cdots$ be such that the corresponding $f^j$ is in the null-space of $X_{g_j}  =X_j$ and $\|F^j\|_{\mathcal{N}^s} = 1$, with $g_j$ within $1/j$ of the Minkowski metric. By the above inequality, $1\leq C_\rho\|F^j\|_{\mathcal{N}^{s-\rho}}.$ Now, $F^j$ has a $\mathcal{N}^s$-weakly convergent subsequence, not shown in notation, to some $F\in \mathcal{N}^s$, which thus strongly converges in $\mathcal{N}^{s-\rho}$. By  the above inequality, $F\neq 0$. But $0 = X_j  f$ converges to $X_M  f$ e.g. in the sense of distributions. So $X_M   f =0$ which by the injectivity of $X_M$, implies that $f = 0$. So we get $F = 0$ a contradiction. 
This shows the injectivity of $X_\delta$ and finishes the proof of Theorem \ref{thm-main3}. 

\epf

 \section*{Acknowledgments}
The authors thank Plamen Stefanov and Gunther Uhlmann for helpful discussions. They also thank the anonymous referees for valuable comments.   A.V. acknowledges support from the National Science Foundation under grants number  DMS-1664683. 


\end{document}